\documentclass[10pt]{amsart}
\usepackage[T1]{fontenc}
\usepackage{lmodern}
\usepackage{geometry}
\usepackage[latin1] {inputenc}
\usepackage{amsmath}
\usepackage{amsfonts, amssymb, textcomp}
\usepackage[colorlinks=flase, linkcolor=red,urlcolor=green, citecolor=blue]{hyperref}
\usepackage{subeqnarray}

\usepackage{latexsym}
\usepackage{fancyhdr}
\usepackage{longtable}
\usepackage{amsmath, amssymb}
\usepackage{graphicx}

\numberwithin{equation}{section}

\theoremstyle{plain}
\newtheorem{theorem}{Theorem}[section]
\newtheorem{lemma}[theorem]{Lemma}
\newtheorem{proposition}[theorem]{Proposition}
\newtheorem{corollary}[theorem]{Corollary}
\theoremstyle{definition}

\newtheorem{remark}[theorem]{Remark}

\newcommand{\RR}{\mathbb{R}}
\newcommand{\CC}{\mathbb{C}}
\newcommand{\NN}{\mathbb{N}}

\newcommand{\supp}{\operatorname{supp}}

\newcommand{\dist}{\operatorname{dist}}

\let\on=\operatorname

\newenvironment{demo}[1]{\par\smallskip\noindent{\bf #1.}}{\par\smallskip}

\title[Sectorial extensions for ultraholomorphic classes defined by weight functions]
{Sectorial extensions for ultraholomorphic classes defined by weight functions}

\author[J.~Jim\'{e}nez-Garrido, J.~Sanz, and G.~Schindl]{Javier Jim\'{e}nez-Garrido, Javier Sanz and Gerhard Schindl}

\begin{document}
\begin{abstract}
We prove an extension theorem for ultraholomorphic classes defined by so-called Braun-Meise-Taylor weight functions $\omega$ and transfer the proofs from the single weight sequence case from V. Thilliez \cite{Thilliezdivision} to the weight function setting. We are following a different approach than the results obtained in \cite{sectorialextensions}, more precisely we are working with real methods by applying the ultradifferentiable Whitney-extension theorem. We are treating both the Roumieu and the Beurling case, the latter one is obtained by a reduction from the Roumieu case.
\end{abstract}

\date{\today}

\maketitle

\section{Introduction}

In a very recent work by the authors~\cite{sectorialextensions}, Roumieu-like ultraholomorphic classes of functions in unbounded sectors of the Riemann surface of the logarithm have been defined by means of so-called Braun-Meise-Taylor weight functions $\omega$ (see \cite{BraunMeiseTaylor90}), and the surjectivity of the Borel map (via the existence of right inverses for this map) has been proved for narrow enough sectors. We refer to the introduction of~\cite{sectorialextensions} for a motivation of this problem and its close links with similar Whitney-like extension problems in the real-variable, ultradifferentiable setting. Despite this connection, no such Whitney extension result was used in the approach followed in~\cite{sectorialextensions}, but instead the techniques employed, of a complex-analytic nature, rested on the use of truncated Laplace transforms whose integral kernel is obtained from optimal flat functions in the corresponding classes. However, the question remained open about the possibility of closely following and generalizing the ideas of V. Thilliez~\cite{Thilliezdivision} in his proof of similar extension results for Denjoy-Carleman ultraholomorphic classes, defined by means of a strongly regular weight sequence $M$, by using different variants of Whitney's theorem.
Our main aim in this paper is to show that Thilliez's techniques are also working in this situation, and how one may also prove the surjectivity of the Borel map for the corresponding Beurling-like classes, which were not considered in~\cite{sectorialextensions}.

In the main result, Theorem \ref{Theorem321}, we prove the surjectivity of the Borel map in ultraholomorphic classes defined by weight functions $\omega$, satisfying several standard assumptions and such that $\gamma(\omega)>1$, and in sectors of opening smaller than $\pi(\gamma(\omega)-1)$. So, as it happens in Thilliez's result, the opening of the sectors for which the result applies is controlled by a growth index $\gamma(\omega)$, which has also been introduced in \cite{sectorialextensions} and is studied in detail in \cite{firstindexpaper}. Moreover, in \cite{firstindexpaper} we consider other indices for $\omega$ and also their relation to the indices $\gamma(M)$ of Thilliez or $\omega(M)$ introduced in \cite{Sanzflatultraholomorphic}.

Some specific points in our arguments deserve mention.
Firstly, we should note that the operations of multiplying or dividing a weight sequence $M=(M_p)_p$ by the factorials, getting $(p!M_p)_p$ or $(M_p/p!)_p$, increase, respectively decrease, the value of the index $\gamma(M)$ by 1, so we may say the first operation takes the sequence to an upper level, while the second one takes it to a lower one. These operations have their counterpart when one considers weight functions $\omega$, expressed by means of the so-called Legendre upper and lower envelopes (or conjugates), and again they increase or decrease the index $\gamma(\omega)$ by 1 (see Section \ref{conugateofaweight}, also \cite[Section 3.1]{sectorialextensions}). It turns out that, depending on the particular setting one considers, it is more natural and/or convenient to express the results departing from the proper level to which the weight function belongs or from the corresponding lower level. Indeed, while in \cite{sectorialextensions} it was preferable to start from the lower level, here we have better considered the upper level as starting point. The upper conjugate  has already played a prominent role in several other extension results, e.g. see \cite{PetzscheVogt}, \cite{BonetBraunMeiseTaylorWhitneyextension}, \cite{whitneyextensionmixedweightfunction}.

Secondly, it is important to highlight that A. Rainer and the third author \cite{dissertation,compositionpaper} have considered ultradifferentiable classes defined by weight matrices (i.e. families of weight sequences), what strictly includes both the Denjoy-Carleman and the Braun-Meise-Taylor approaches. Since with each weight function we can associate a canonical weight matrix $\Omega=\{W^x: x>0\}$ defining the same corresponding classes, the idea is to transfer proofs from the weight sequence setting to the weight function one by working with the sequences $W^x$, and this will be fruitful in our regards.
A main challenge in this respect is that one can not use directly the proofs from \cite{Thilliezdivision} and replace there the weight sequence $M$ by some/each $W^x\in\Omega$. The reason is that, in general, we can not assume some/each $W^x$ to be strongly regular as needed in \cite{Thilliezdivision}, since this would lead us to the uninteresting case of a constant weight matrix,
see the end of Section \ref{ultraholomorphicroumieu} for precise explanations. Moreover, in case the sector is not a subset of $\CC$, i.e. if its opening is greater than $2\pi$, we have to consider a special ramification construction for (associated) weight matrices, discussed in Section \ref{ramificationpreparation}, and which is intimately related to the results and techniques from \cite[Section 5]{sectorialextensions}.

\vspace{6pt}

The paper is organized as follows. In Section \ref{basic} we are collecting all the preliminary, mostly well-known, information and notation concerning weight sequences, weight functions and weight matrices, and we are defining the ultraholomorphic classes under consideration. In Section \ref{associatedweightsandconjugates} we are introducing the notion of weight functions associated to weight sequences and recalling some basic facts about Legendre (also called Young) conjugates of functions. Several auxiliary results are stated, most of them have already been shown in \cite[Sections 2, 3]{sectorialextensions}.

Section \ref{growthindexgamma} is devoted to the definition of the growth index $\gamma(\omega)$ in terms of a weight function and we state there the results which are playing the key role in the forthcoming sections. As mentioned before, precise proofs, the connection to $\gamma(M)$ and much more information on this topic is given in \cite{firstindexpaper}. In Section \ref{sectoriallyflat} we construct the optimal flat functions, for which some preparatory work has already been done in \cite[Section 6]{sectorialextensions}.

In Section \ref{ultraholomorphicextension}, we state and prove the central result Theorem \ref{Theorem321}.
Finally, in Section \ref{Beurlingcase} we succeed to treat the ultraholomorphic Beurling-like classes, which have not been considered in \cite{sectorialextensions} and, even in the weight sequence setting, have not been used very frequently in the literature (e.g. see \cite[Section 3.4]{Thilliezdivision} and \cite{Schmetsvaldivia00}). We reduce the proof of the corresponding extension result to the Roumieu case, a method which has already been used several times, e.g. see \cite[Section 3.4]{Thilliezdivision}, \cite[Lemma 4.4, Thm. 4.5]{BonetBraunMeiseTaylorWhitneyextension}).

In the appendix we are investigating the special class of weight functions satisfying a quite strong condition, denoted by \hyperlink{om7}{$(\omega_7)$} in this work and that has already appeared, in an equivalent form, in the work of U. Franken~\cite{Franken95}. Such weights always fulfill $\gamma(\omega)=+\infty$, hence the main result is valid for sectors of any opening. Moreover, the technical ramification construction in Section \ref{ramificationpreparation} becomes superfluous in this case, due to the strong properties verified in this case by the weight matrix associated with $\omega$.

\section{Basic definitions}\label{basic}
We will use the following general notation: $\mathcal{E}$ shall denote the class of all smooth functions and $\mathcal{H}$ the class of holomorphic functions. We will write $\NN_{>0}=\{1,2,\dots\}$ and $\NN=\NN_{>0}\cup\{0\}$, moreover let $\RR_{>0}:=\{x\in\RR: x>0\}$ denote the set of all positive real numbers. For a given multiindex $k=(k_1,\dots,k_n)\in\NN^n$ we put $|k|:=k_1+\dots+k_n$.

\subsection{Weight sequences}
A sequence $M=(M_k)_k\in\RR_{>0}^{\NN}$ is called a {\itshape weight sequence}, define also $m=(m_k)_k$ by $m_k:=\frac{M_k}{k!}$ and $\mu=(\mu_k)_k$ by $\mu_k:=\frac{M_k}{M_{k-1}}$, $\mu_0:=1$.

$M$ is called {\itshape normalized} if $1=M_0\le M_1$ and which can always be assumed without loss of generality.

$(1)$ $M$ is {\itshape log-convex}, if
$$\hypertarget{lc}{(\text{lc})}:\Leftrightarrow\;\forall\;j\in\NN_{>0}:\;M_j^2\le M_{j-1} M_{j+1}$$
and {\itshape strongly log-convex}, if
$$\hypertarget{slc}{(\text{slc})}:\Leftrightarrow\;\forall\;j\in\NN_{>0}:\;m_j^2\le m_{j-1} m_{j+1}.$$
If $M$ is log-convex and normalized, then $M$ and $k\mapsto(M_k)^{1/k}$ are both nondecreasing and $M_jM_k\le M_{j+k}$ holds for all $j,k\in\NN$, e.g. see \cite[Remark 2.0.3, Lemmata 2.0.4, 2.0.6]{diploma}.

$(2)$ $M$ has {\itshape moderate growth} if
$$\hypertarget{mg}{(\text{mg})}:\Leftrightarrow\exists\;C\ge 1\;\forall\;j,k\in\NN:\;M_{j+k}\le C^{j+k}M_j M_k.$$
Given two (different) weight sequences we will also consider the following mixed moderate growth condition
$$\hypertarget{genmg}{(M,N)_{(\on{mg})}}:\Leftrightarrow\;\exists\;C\ge 1\;\forall\;j,k\in\NN:\;\;\; M_{j+k}\le C^{j+k}N_jN_k.$$
Note that $(M,N)_{(\on{mg})}\Leftrightarrow(m,n)_{(\on{mg})}$ by changing the constant $C$.

$(3)$ $M$ is called {\itshape nonquasianalytic,} if
$$\hypertarget{mnq}{(\text{nq})}:\Leftrightarrow\;\sum_{p=1}^{+\infty}\frac{M_{p-1}}{M_p}<+\infty.$$
If $M$ is log-convex, then using {\itshape Carleman's inequality} one can show that $\sum_{p=1}^{+\infty}\frac{M_{p-1}}{M_p}<+\infty\Leftrightarrow\sum_{p=1}^{+\infty}\frac{1}{(M_p)^{1/p}}<+\infty$ (for a proof e.g. see \cite[Proposition 4.1.7]{diploma} and \cite[Theorem 1.3.8]{hoermander}).

$(4)$ $M$ has $(\beta_1)$ if
$$\hypertarget{beta1}{(\beta_1)}:\Leftrightarrow\;\exists\;Q\in\NN_{>0}:\;\liminf_{p\rightarrow+\infty}\frac{\mu_{Qp}}{\mu_p}>Q,$$
and $(\gamma_1)$ if
$$\hypertarget{gamma1}{(\gamma_1)}:\Leftrightarrow\sup_{p\in\NN_{>0}}\frac{\mu_p}{p}\sum_{k\ge p}\frac{1}{\mu_k}<+\infty.$$
By \cite[Proposition 1.1]{petzsche} both conditions are equivalent for log-convex $M$ and for this proof condition \hyperlink{mnq}{$(\text{nq})$}, which is a general assumption in \cite{petzsche}, was not necessary. In the literature a sequence satisfying \hyperlink{gamma1}{$(\gamma_1)$} is also called ``strongly nonquasianalytic''.

Due to technical reasons it is often convenient to assume several properties for $M$ at the same time and hence we define the class\vspace{6pt}

\centerline{$M\in\hypertarget{SRset}{\mathcal{SR}}$, if $M$ is normalized and has \hyperlink{slc}{$(\text{slc})$}, \hyperlink{mg}{$(\text{mg})$} and \hyperlink{gamma1}{$(\gamma_1)$}.}\vspace{6pt}

Using this notation we see that $M\in\hyperlink{SRset}{\mathcal{SR}}$ if and only if $m$ is a {\itshape strongly regular sequence} in the sense of \cite[1.1]{Thilliezdivision} (and this terminology has also been used by several authors so far, e.g. see \cite{Sanzflatultraholomorphic}, \cite{Sanzsummability}). At this point we want to make the reader aware that here we are using the same notation as it has already been used by the authors in \cite{sectorialextensions}, whereas in \cite{Thilliezdivision} and also in \cite{firstindexpaper} the sequence $M$ is precisely $m$ in the notation in this work. We prefer to work with $M$ instead of $m$ since we want to use the sequences contained in the associated weight matrix with a given weight function, see \ref{weightmatrixfromfunction} below for more explanations.


We write $M\le N$ if and only if $M_p\le N_p$ holds for all $p\in\NN$ and define
$$M\hypertarget{precsim}{\precsim}N:\Leftrightarrow\;\exists\;C\ge 1\;\forall\;p\in\NN:\; M_p\le C^pN_p\Longleftrightarrow\sup_{p\in\NN_{>0}}\left(\frac{M_p}{N_p}\right)^{1/p}<+\infty.$$
$M$ and $N$ are called {\itshape equivalent} if
$$M\hypertarget{approx}{\approx}N:\Leftrightarrow\;M\hyperlink{precsim}{\precsim}N\;\text{and}\;N\hyperlink{precsim}{\precsim}M.$$
Note that $M\precsim N\Leftrightarrow m\precsim n$ by changing the constant.

For convenience we define the set
$$\hypertarget{LCset}{\mathcal{LC}}:=\{M\in\RR_{>0}^{\NN}:\;M\;\text{normalized, log-convex},\;\lim_{k\rightarrow+\infty}(M_k)^{1/k}=+\infty\}.$$

\subsection{Weight functions $\omega$ in the sense of Braun-Meise-Taylor}\label{weightfunctionclasses}
A function $\omega:[0,+\infty)\rightarrow[0,+\infty)$ is called a {\itshape weight function} (in the terminology of \cite[Section 2.2]{sectorialextensions}), if it is continuous, nondecreasing, $\omega(0)=0$ and $\lim_{t\rightarrow+\infty}\omega(t)=+\infty$. If $\omega$ satisfies in addition $\omega(t)=0$ for all $t\in[0,1]$, then we call $\omega$ a {\itshape normalized weight function} and for convenience we will write that $\omega$ has $\hypertarget{om0}{(\omega_0)}$ if it is a normalized weight.

Given an arbitrary function $\omega$ we will denote by $\omega^{\iota}$ the function $\omega^{\iota}(t):=\omega(1/t)$ for any $t>0$.\vspace{6pt}

Moreover we consider the following conditions, this list of properties has already been used in ~\cite{dissertation}.


\begin{itemize}
\item[\hypertarget{om1}{$(\omega_1)}$] $\omega(2t)=O(\omega(t))$ as $t\rightarrow+\infty$, i.e. $\exists\;L\ge 1\;\forall\;t\ge 0:\;\;\;\omega(2t)\le L(\omega(t)+1)$.

\item[\hypertarget{om2}{$(\omega_2)$}] $\omega(t)=O(t)$ as $t\rightarrow+\infty$.

\item[\hypertarget{om3}{$(\omega_3)$}] $\log(t)=o(\omega(t))$ as $t\rightarrow+\infty$ ($\Leftrightarrow\lim_{t\rightarrow+\infty}\frac{t}{\varphi_{\omega}(t)}=0$).

\item[\hypertarget{om4}{$(\omega_4)$}] $\varphi_{\omega}:t\mapsto\omega(e^t)$ is a convex function on $\RR$.

\item[\hypertarget{om5}{$(\omega_5)$}] $\omega(t)=o(t)$ as $t\rightarrow+\infty$.

\item[\hypertarget{om6}{$(\omega_6)$}] $\exists\;H\ge 1\;\forall\;t\ge 0:\;2\omega(t)\le\omega(H t)+H$.

\item[\hypertarget{om7}{$(\omega_7)$}] $\exists\;H>0\;\exists\;C>0\;\forall\;t\ge 0: \omega(t^2)\le C\omega(H t)+C.$

\item[\hypertarget{omnq}{$(\omega_{\text{nq}})$}] $\int_1^{+\infty}\frac{\omega(t)}{t^2}dt<+\infty.$

\item[\hypertarget{omsnq}{$(\omega_{\text{snq}})$}] $\exists\;C>0:\;\forall\;y>0: \int_1^{+\infty}\frac{\omega(y t)}{t^2}dt\le C\omega(y)+C$.
\end{itemize}

An interesting example is $\sigma_s(t):=\max\{0,\log(t)^s\}$, $s>1$, which satisfies all listed properties except \hyperlink{om6}{$(\omega_6)$}. It is well-known that the ultradifferentiable class defined by using the weight $t\mapsto t^{1/s}$ yields the Gevrey class $G^s=(p!^s)_{p\in\NN}$ of index $s>1$ (see Section \ref{ultradifferentiableclasses} below for precise definitions), it satisfies all listed properties except \hyperlink{om7}{$(\omega_7)$}.

Note that for any weight function property \hyperlink{omnq}{$(\omega_{\text{nq}})$} implies \hyperlink{om5}{$(\omega_5)$} because $\int_t^{+\infty}\frac{\omega(u)}{u^2}du\ge\omega(t)\int_t^{+\infty}\frac{1}{u^2}du=\frac{\omega(t)}{t}$.\vspace{6pt}

For convenience we define the sets
$$\hypertarget{omset0}{\mathcal{W}_0}:=\{\omega:[0,\infty)\rightarrow[0,\infty): \omega\;\text{has}\;\hyperlink{om0}{(\omega_0)},\hyperlink{om3}{(\omega_3)},\hyperlink{om4}{(\omega_4)}\},\hspace{20pt}\hypertarget{omset1}{\mathcal{W}}:=\{\omega\in\mathcal{W}_0: \omega\;\text{has}\;\hyperlink{om1}{(\omega_1)}\}.$$

For any $\omega\in\hyperlink{omset0}{\mathcal{W}_0}$ we define the {\itshape Legendre-Fenchel-Young-conjugate} of $\varphi_{\omega}$ by
\begin{equation}\label{legendreconjugate}
\varphi^{*}_{\omega}(x):=\sup\{x y-\varphi_{\omega}(y): y\ge 0\},\;\;\;x\ge 0,
\end{equation}
with the following properties, e.g. see \cite[Remark 1.3, Lemma 1.5]{BraunMeiseTaylor90}: It is convex and nondecreasing, $\varphi^{*}_{\omega}(0)=0$, $\varphi^{**}_{\omega}=\varphi_{\omega}$, $\lim_{x\rightarrow+\infty}\frac{x}{\varphi^{*}_{\omega}(x)}=0$ and finally $x\mapsto\frac{\varphi_{\omega}(x)}{x}$ and $x\mapsto\frac{\varphi^{*}_{\omega}(x)}{x}$ are nondecreasing on $[0,+\infty)$. For any $\omega\in\hyperlink{omset1}{\mathcal{W}}$ we can assume w.l.o.g. that $\omega$ is $\mathcal{C}^1$ (see \cite[Lemma 1.7]{BraunMeiseTaylor90}).\vspace{6pt}

Let $\sigma,\tau$ be weight functions, we write
$$\sigma\hypertarget{ompreceq}{\preceq}\tau:\Leftrightarrow\tau(t)=O(\sigma(t))\;\text{as}\;t\rightarrow+\infty$$
and call them equivalent if
$$\sigma\hypertarget{sim}{\sim}\tau:\Leftrightarrow\sigma\hyperlink{ompreceq}{\preceq}\tau\;\text{and}\;\tau\hyperlink{ompreceq}{\preceq}\sigma.$$
Motivated by the notion of a strong weight function given in \cite{BonetBraunMeiseTaylorWhitneyextension}\vspace{6pt}

\centerline{$\omega$ will be called a strong weight, if $\omega\in\hyperlink{omset0}{\mathcal{W}_0}$ and in addition \hyperlink{omsnq}{$(\omega_{\text{snq}})$} is satisfied.}\vspace{6pt}

We recall \cite[Proposition 1.3]{MeiseTaylor88}, where \hyperlink{omsnq}{$(\omega_{\text{snq}})$} was characterized, and \cite[Corollary 1.4]{MeiseTaylor88}:

\begin{proposition}\label{Prop13MT88}
Let $\omega:[0,+\infty)\longrightarrow[0,+\infty)$ be a weight function, then the following are equivalent:
\begin{itemize}
\item[$(i)$] $\lim_{\varepsilon\rightarrow 0}\limsup_{t\rightarrow+\infty}\frac{\varepsilon\omega(t)}{\omega(\varepsilon t)}=0$,
\item[$(ii)$] $\exists\;K>1$ such that $\limsup_{t\rightarrow+\infty}\frac{\omega(Kt)}{\omega(t)}<K$,
\item[$(iii)$] $\omega$ satisfies \hyperlink{omsnq}{$(\omega_{\on{snq}})$},
\item[$(iv)$] there exists a nondecreasing concave function $\kappa: [0,+\infty)\longrightarrow[0,+\infty)$ such that $\omega\hyperlink{sim}{\sim}\kappa$ and $\kappa$ satisfies \hyperlink{omsnq}{$(\omega_{\on{snq}})$}. More precisely $\kappa=\kappa_{\omega}$ with
    \begin{equation*}\label{kappa}
\kappa_{\omega}(t):=\int_1^{+\infty}\frac{\omega(tu)}{u^2}du=t\int_t^{+\infty}\frac{\omega(u)}{u^2}du,\;\;\;\forall\;t>0,\hspace{20pt}\kappa_{\omega}(0):=0.
\end{equation*}
\end{itemize}
Consequently $\omega$ has also \hyperlink{om1}{$(\omega_1)$} and \hyperlink{om5}{$(\omega_5)$}. If $\omega$ satisfies one of the equivalent conditions above, then there exists some $0<\alpha<1$ such that $\omega(t)=O(t^{\alpha})$ as $t\rightarrow+\infty$.

It is well-known that all properties from \hyperlink{omset1}{$\mathcal{W}$} except normalization can be transferred from $\omega$ to $\kappa_{\omega}$, e.g. see \cite[Remark 3.2]{BonetMeiseTaylorSurjectivity}.
\end{proposition}

Note that concavity and $\omega(0)=0$ imply subadditivity and this yields \hyperlink{om1}{$(\omega_1)$}.

W.l.o.g. in $(ii)$ above we can assume $K>p$ for any $p\in\NN_{\ge 2}$ because if $1<K<2$ then iterating $(ii)$ yields
\begin{equation}\label{omegasnqiteration}
\limsup_{t\rightarrow+\infty}\frac{\omega(K^2t)}{\omega(t)}=\limsup_{t\rightarrow+\infty}\frac{\omega(K^2t)}{\omega(Kt)}\frac{\omega(Kt)}{\omega(t)}\le\limsup_{t\rightarrow+\infty}\frac{\omega(K^2t)}{\omega(Kt)}\limsup_{t\rightarrow+\infty}\frac{\omega(Kt)}{\omega(t)}<K^2,
\end{equation}
and so after $n(p)$ steps we get $K^n>p$.\vspace{6pt}

As in \cite{sectorialextensions} for any given weight function $\omega$ and $s>0$ we define a new weight function $\omega^s$ by
\begin{equation}\label{defiRamifiedWeightFunction}
\omega^s(t):=\omega(t^s),\quad t\ge 0.
\end{equation}
If $\omega$ satisfies any of the properties \hyperlink{om1}{$(\omega_1)$}, \hyperlink{om3}{$(\omega_3)$}, \hyperlink{om4}{$(\omega_4)$}, \hyperlink{om6}{$(\omega_6)$} or \hyperlink{om7}{$(\omega_7)$}, then the same holds for $\omega^s$ (and also normalization is preserved). But \hyperlink{om5}{$(\omega_5)$} and \hyperlink{omsnq}{$(\omega_{\text{snq}})$} might not be preserved. Indeed, this last fact motivates the introduction of the index $\gamma(\omega)$ in Section \ref{growthindexgamma} below.

\subsection{Weight matrices}\label{classesweightmatrices}
For the following definitions and conditions see also \cite[Section 4]{compositionpaper}.

Let $\mathcal{I}=\RR_{>0}$ denote the index set (equipped with the natural order), a {\itshape weight matrix} $\mathcal{M}$ associated with $\mathcal{I}$ is a (one parameter) family of weight sequences $\mathcal{M}:=\{M^x\in\RR_{>0}^{\NN}: x\in\mathcal{I}\}$, such that
$$\hypertarget{Marb}{(\mathcal{M})}:\Leftrightarrow\;\forall\;x\in\mathcal{I}:\;M^x\;\text{is normalized, nondecreasing},\;M^{x}\le M^{y}\;\text{for}\;x\le y.$$
We call a weight matrix $\mathcal{M}$ {\itshape standard log-convex,} if
$$\hypertarget{Msc}{(\mathcal{M}_{\on{sc}})}:\Leftrightarrow(\mathcal{M})\; \text{and}\;\forall\;x\in\mathcal{I}:\;M^x\in\hyperlink{LCset}{\mathcal{LC}}.$$
Moreover, we put $m^x_p:=\frac{M^x_p}{p!}$ for $p\in\NN$, and $\mu^x_p:=\frac{M^x_p}{M^x_{p-1}}$ for $p\in\NN_{>0}$, $\mu^x_0:=1$.

A matrix is called {\itshape constant} if $\mathcal{M}=\{M\}$ or more generally if $M^x\hyperlink{approx}{\approx}M^y$ for all $x,y\in\mathcal{I}$.

We are going to consider the following properties for $\mathcal{M}$ (see \cite[Section 4.1]{compositionpaper} and \cite[Section 7.2]{dissertation}):\par\vskip.3cm

\hypertarget{R-mg}{$(\mathcal{M}_{\{\text{mg}\}})$} \hskip1cm $\forall\;x\in\mathcal{I}\;\exists\;C>0\;\exists\;y\in\mathcal{I}\;\forall\;j,k\in\NN: M^x_{j+k}\le C^{j+k} M^y_j M^y_k$.\par\vskip.3cm

\hypertarget{B-mg}{$(\mathcal{M}_{(\text{mg})})$} \hskip1cm $\forall\;x\in\mathcal{I}\;\exists\;C>0\;\exists\;y\in\mathcal{I}\;\forall\;j,k\in\NN: M^y_{j+k}\le C^{j+k} M^x_j M^x_k$.\par\vskip.3cm
\hypertarget{R-L}{$(\mathcal{M}_{\{\text{L}\}})$} \hskip1cm $\forall\;C>0\;\forall\;x\in\mathcal{I}\;\exists\;D>0\;\exists\;y\in\mathcal{I}\;\forall\;k\in\NN: C^k M^x_k\le D M^y_k$.\vskip.3cm
\hypertarget{B-L}{$(\mathcal{M}_{(\text{L})})$} \hskip1cm $\forall\;C>0\;\forall\;x\in\mathcal{I}\;\exists\;D>0\;\exists\;y\in\mathcal{I}\;\forall\;k\in\NN: C^k M^y_k\le D M^x_k$.\vskip.3cm

Let $\mathcal{M}=\{M^x: x\in\mathcal{I}\}$ and $\mathcal{N}=\{N^x: x\in\mathcal{J}\}$ be \hyperlink{Marb}{$(\mathcal{M})$}, define
$$\mathcal{M}\hypertarget{Mroumprecsim}{\{\precsim\}}\mathcal{N}\;:
\Leftrightarrow\forall\;x\in\mathcal{I}\;\exists\;y\in\mathcal{J}:\;
M^x\hyperlink{precsim}{\precsim}N^y,$$
and call them equivalent if
$$\mathcal{M}\hypertarget{Mroumapprox}{\{\approx\}}\mathcal{N}\;:
\Leftrightarrow\mathcal{M}\hyperlink{Mroumprecsim}{\{\precsim\}}\mathcal{N}\;\text{and}\;\mathcal{N}\hyperlink{Mroumprecsim}{\{\precsim\}}\mathcal{M}.$$

\subsection{Weight matrices obtained by weight functions}\label{weightmatrixfromfunction}
We summarize some facts which are shown in \cite[Section 5]{compositionpaper} and will be needed in this work.

\begin{itemize}
\item[$(i)$] The idea was that to each $\omega\in\hyperlink{omset1}{\mathcal{W}}$ we can associate a \hyperlink{Msc}{$(\mathcal{M}_{\text{sc}})$} weight matrix $\Omega:=\{W^l=(W^l_j)_{j\in\NN}: l>0\}$ by\vspace{6pt}

    \centerline{$W^l_j:=\exp\left(\frac{1}{l}\varphi^{*}_{\omega}(lj)\right)$.}\vspace{6pt}

In general it is not clear that $W^x$ is strongly log-convex, i.e. $w^x$ is log-convex, too.

\item[$(ii)$] $\Omega$ satisfies \hyperlink{R-mg}{$(\mathcal{M}_{\{\text{mg}\}})$} and \hyperlink{B-mg}{$(\mathcal{M}_{(\text{mg})})$}, more precisely
    \begin{equation}\label{newmoderategrowth}
    \forall\;l>0\;\forall\;j,k\in\NN:\;\;\;W^l_{j+k}\le W^{2l}_jW^{2l}_k,
    \end{equation}
    i.e. \hyperlink{genmg}{$(W^l,W^{2l})_{(\on{mg})}$}. $\Omega$ has also both \hyperlink{R-L}{$(\mathcal{M}_{\{\text{L}\}})$} and \hyperlink{B-L}{$(\mathcal{M}_{(\text{L})})$}, more precisely we get
     \begin{equation}\label{newexpabsorb}
     \forall\;h\ge 1\;\exists\;A\ge 1\;\forall\;l>0\;\exists\;D\ge 1\;\forall\;j\in\NN:\;\;\;h^jW^l_j\le D W^{Al}_j.
     \end{equation}
     In fact we can take $A=L_1^a$, where $L_1=L(L+1)\ge L$ and $L$ is the constant arising in \hyperlink{om1}{$(\omega_1)$}, the number $a\in\NN_{>0}$ is chosen minimal to have $\exp(a)\ge h$ (see the proof of \cite[Proposition 3.3.1, p.18]{dissertation} and of \cite[Lemma 5.9 $(5.10)$]{compositionpaper}).

\item[$(iii)$] Equivalent weight functions yield equivalent associated weight matrices. Note that \hyperlink{R-mg}{$(\mathcal{M}_{\{\text{mg}\}})$} is stable with respect to \hyperlink{Mroumapprox}{$\{\approx\}$}, whereas \hyperlink{R-L}{$(\mathcal{M}_{\{\text{L}\}})$} not necessarily.

\item[$(iv)$] \hyperlink{om5}{$(\omega_5)$} holds if and only if $\lim_{p\rightarrow+\infty}(w^l_p)^{1/p}=+\infty$ for all $l>0$.

\end{itemize}
In fact all statements from above except \eqref{newexpabsorb} are even true for $\omega\in\hyperlink{omset0}{\mathcal{W}_0}$.

Moreover we have:
\begin{remark}\label{importantremark}
Let $\omega\in\hyperlink{omset0}{\mathcal{W}_0}$ be given, then $\omega$ satisfies
 \begin{itemize}
 \item[$(a)$] \hyperlink{omnq}{$(\omega_{\on{nq}})$} if and only if some/each $W^l$ satisfies \hyperlink{mnq}{$(\on{nq})$},
 \item[$(b)$] \hyperlink{om6}{$(\omega_6)$} if and only if some/each $W^l$ satisfies \hyperlink{mg}{$(\on{mg})$} if and only if $W^l\hyperlink{approx}{\approx}W^n$ for each $l,n>0$. Consequently \hyperlink{om6}{$(\omega_6)$} is characterizing the situation when $\Omega$ is constant, i.e. $\Omega=\{W^l\}$ for some/each $l>0$.
 \end{itemize}
\end{remark}

Finally the next calculation relates the matrices associated to $\omega$ and $\omega^s$, $s>0$, defined in \eqref{defiRamifiedWeightFunction}. Since $\varphi_{\omega^s}(t)=\varphi_{\omega}(st)$ we get $\varphi^{*}_{\omega^s}(x)=\varphi^{*}_{\omega}(x/s)$ for any $x\ge 0$ and $s>0$. Hence by denoting $V^{l,s}$ the $l$-th sequence of the weight matrix associated to $\omega^s$ we have for all $j\in\NN$:
\begin{equation*}\label{weightmatrixV}
V^{l,s}_j=\exp\left(\frac{1}{l}\varphi^{*}_{\omega^s}(lj)\right)=\exp\left(\frac{1}{l}\varphi^{*}_{\omega}(lj/s)\right)=\exp\left(\frac{s}{ls}\varphi^{*}_{\omega}(lj/s)\right)=\left(W^{l/s}_j\right)^{1/s}.
\end{equation*}

\subsection{Classes of ultraholomorphic functions}\label{ultraholomorphicroumieu}
We introduce now the classes under consideration in this paper, see also \cite[Section 2.5]{sectorialextensions}. For the following definitions, notation and more details we refer to \cite[Section 2]{Sanzflatultraholomorphic}. Let $\mathcal{R}$ be the Riemann surface of the logarithm. We wish to work in general unbounded sectors in $\mathcal{R}$ with vertex at $0$, but all our results will be unchanged under rotation, so we will only consider sectors bisected by direction $0$: For $\gamma>0$ we set $$S_{\gamma}:=\left\{z\in\mathcal{R}: |\arg(z)|<\frac{\gamma\pi}{2}\right\},$$
i.e. the unbounded sector of opening $\gamma\pi$, bisected by direction $0$.

Let $M$ be a weight sequence, $S\subseteq\mathcal{R}$ an (unbounded) sector and $h>0$. We define
$$\mathcal{A}_{M,h}(S):=\{f\in\mathcal{H}(S): \|f\|_{M,h}:=\sup_{z\in S, p\in\NN}\frac{|f^{(p)}(z)|}{h^p M_p}<+\infty\}.$$
$(\mathcal{A}_{M,h}(S),\|\cdot\|_{M,h})$ is a Banach space and we put
\begin{equation*}
\mathcal{A}_{\{M\}}(S):=\bigcup_{h>0}\mathcal{A}_{M,h}(S).
\end{equation*}
$\mathcal{A}_{\{M\}}(S)$ is called the Denjoy-Carleman ultraholomorphic class (of Roumieu type) associated with $M$ in the sector $S$ (it is a $(LB)$ space). Analogously we introduce the space of complex sequences
$$\Lambda^n_{M,h}:=\{a=(a_p)_p\in\CC^{\NN^n}: |a|_{M,h}:=\sup_{p\in\NN^n}\frac{|a_p|}{h^{|p|} M_{|p|}}<+\infty\}$$
and put $\Lambda^n_{\{M\}}:=\bigcup_{h>0}\Lambda^n_{M,h}$. The (asymptotic) {\itshape Borel map} $\mathcal{B}$ is given by
\begin{equation*}
\mathcal{B}:\mathcal{A}_{\{M\}}(S)\longrightarrow\Lambda^1_{\{M\}},\hspace{15pt}f\mapsto(f^{(p)}(0))_{p\in\NN},
\end{equation*}
where $f^{(p)}(0):=\lim_{z\in S, z\rightarrow 0}f^{(p)}(z)$.\vspace{6pt}

Not very frequently in the literature also Denjoy-Carleman ultraholomorphic classes $\mathcal{A}_{(M)}(S)$ (and analogously sequence classes $\Lambda^n_{(M)}$) of Beurling type have been investigated, e.g. see \cite[Section 3.4]{Thilliezdivision} and \cite{Schmetsvaldivia00}. In this case we put
$$\mathcal{A}_{(M)}(S):=\bigcap_{h>0}\mathcal{A}_{M,h}(S),\hspace{30pt}\Lambda^n_{(M)}:=\bigcap_{h>0}\Lambda^n_{M,h},$$
so consider the projective structure with respect to $h>0$ and obtain {\itshape Fr\'echet spaces}.\vspace{6pt}

Similarly as for the ultradifferentiable case, we now define ultraholomorphic classes associated with $\omega\in\hyperlink{omset0}{\mathcal{W}_0}$. Given an unbounded sector $S$, and for every $l>0$, we first define
$$\mathcal{A}_{\omega,l}(S):=\{f\in\mathcal{H}(S): \|f\|_{\omega,l}:=\sup_{z\in S, p\in\NN}\frac{|f^{(p)}(z)|}{\exp(\frac{1}{l}\varphi^{*}_{\omega}(lp))}<+\infty\}.$$
$(\mathcal{A}_{\omega,l}(S),\|\cdot\|_{\omega,l})$ is a Banach space and we put
\begin{equation*}
\mathcal{A}_{\{\omega\}}(S):=\bigcup_{l>0}\mathcal{A}_{\omega,l}(S).
\end{equation*}
$\mathcal{A}_{\{\omega\}}(S)$ is called the Denjoy-Carleman ultraholomorphic class (of Roumieu type) associated with $\omega$ in the sector $S$ (it is a $(LB)$ space). Correspondingly, we introduce the space of complex sequences
$$\Lambda^n_{\omega,l}:=\{a=(a_p)_p\in\CC^{\NN^n}: |a|_{\omega,l}:=\sup_{p\in\NN^n}\frac{|a_p|}{\exp(\frac{1}{l}\varphi^{*}_{\omega}(l|p|))}<+\infty\}$$
and put $\Lambda^n_{\{\omega\}}:=\bigcup_{l>0}\Lambda^n_{\omega,l}$. So in this case we get the {\itshape Borel map} $\mathcal{B}:\mathcal{A}_{\{\omega\}}(S)\longrightarrow\Lambda^1_{\{\omega\}}$.\vspace{6pt}

Again the corresponding Beurling type class $\mathcal{A}_{(\omega)}(S)$ (see Section \ref{Beurlingcase}) is the projective limit of $\mathcal{A}_{\omega,l}(S)$ with respect to $l>0$, similarly for $\Lambda^n_{(\omega)}$ too and so $\mathcal{B}:\mathcal{A}_{(\omega)}(S)\longrightarrow\Lambda^1_{(\omega)}$.\vspace{6pt}

Finally, we recall the ultradifferentiable function classes of {\itshape Roumieu type} defined by a weight matrix $\mathcal{M}$, introduced in \cite[Section 7]{dissertation} and also in \cite[Section 4.2]{compositionpaper}.

Given a weight matrix $\mathcal{M}=\{M^x\in\RR_{>0}^{\NN}: x\in\RR_{>0}\}$ and a sector $S$ we may define ultraholomorphic classes  $\mathcal{A}_{\{\mathcal{M}\}}(S)$  of {\itshape Roumieu type} as
\begin{equation*}
\mathcal{A}_{\{\mathcal{M}\}}(S):=\bigcup_{x\in\RR_{>0}}\mathcal{A}_{\{M^x\}}(S),
\end{equation*}
and accordingly, $\Lambda^n_{\{\mathcal{M}\}}:=\bigcup_{x\in\RR_{>0}}\Lambda^n_{\{M^x\}}$.

Let now $\omega\in\hyperlink{omset1}{\mathcal{W}}$ be given and let $\Omega$ be the associated weight matrix defined in Section \ref{weightmatrixfromfunction} $(i)$, then
\begin{equation}\label{equaEqualitySpacesWeightFunctionMatrix}
\mathcal{A}_{\{\omega\}}(S)=\mathcal{A}_{\{\Omega\}}(S)
\end{equation}
holds as locally convex vector spaces. This equality is an easy consequence of \cite[Lemma 5.9 (5.10)]{compositionpaper} (property \hyperlink{R-L}{$(\mathcal{M}_{\{\text{L}\}})$}) and the way how the seminorms are defined in these spaces. As one also has $\Lambda^n_{\{\omega\}}=\Lambda^n_{\{\Omega\}}$, the Borel map $\mathcal{B}$ makes sense in these last classes,
$\mathcal{B}:\mathcal{A}_{\{\Omega\}}(S)\longrightarrow\Lambda^1_{\{\Omega\}}$.\vspace{6pt}

Analogously as in \cite{dissertation} and \cite{compositionpaper} we can also consider ultraholomorphic classes $\mathcal{A}_{(\mathcal{M})}(S)$ of {\itshape Beurling type} as
\begin{equation*}
\mathcal{A}_{(\mathcal{M})}(S):=\bigcap_{x\in\RR_{>0}}\mathcal{A}_{(M^x)}(S),
\end{equation*}
and accordingly, $\Lambda^n_{(\mathcal{M})}:=\bigcap_{x\in\RR_{>0}}\Lambda^n_{(M^x)}$. Again for any $\omega\in\hyperlink{omset1}{\mathcal{W}}$ we have $\mathcal{A}_{(\omega)}(S)=\mathcal{A}_{(\Omega)}(S)$ as locally convex vector spaces by analogous reasons as in the Roumieu case before and so $\mathcal{B}:\mathcal{A}_{(\Omega)}(S)\longrightarrow\Lambda^1_{(\Omega)}$.\vspace{6pt}

In any of the considered ultraholomorphic classes, an element $f$ is said to be flat if $f^{(p)}(0) = 0$ for
every $p\in\NN$, that is, $\mathcal{B}(f)$ is the null sequence.\vspace{6pt}

With all this information we may claim: In general it seems reasonable to transfer well-known results and proofs from the weight sequence to the weight function setting by using its associated weight matrix and even to the setting in which an abstract weight matrix is given from the start. Unfortunately in many cases it is often impossible to replace $M$ by $W^x$ in the proofs directly since, due to technical reasons, undesirable, respectively too strong, conditions have been imposed on the weight sequences.\vspace{6pt}

More precisely, in \cite{Thilliezdivision}, the general framework has been working with {\itshape strongly regular} weight sequences.

The statements from Section \ref{weightmatrixfromfunction} show that we can not assume in general that some/each $W^x\in\hyperlink{SRset}{\mathcal{SR}}$ (respectively equivalently that some/each $w^x$ is {\itshape strongly regular}). First log-convexity for $w^x$ is not clear in general, second \hyperlink{mg}{$(\on{mg})$} for some/each $W^x$ would yield the constant case. Finally \hyperlink{gamma1}{$(\gamma_1)$} for some/each $W^x$ is also a strong assumption, a mixed condition seems to be more reasonable in the nonconstant case, e.g. see \cite[Thm. 5.12, Cor. 5.13]{whitneyextensionweightmatrix} and \cite[Thm 5.15]{whitneyextensionmixedweightfunction}. In any case it is natural that \hyperlink{omsnq}{$(\omega_{\text{snq}})$} should be the correct assumption and be assumed in the weight function setting.

\section{Associated weight functions and conjugates}\label{associatedweightsandconjugates}
\subsection{Functions $\omega_M$ and $h_M$}\label{assofunction}
Let $M\in\RR_{>0}^{\NN}$ (with $M_0=1$), then the {\itshape associated function} $\omega_M: \RR_{\ge 0}\rightarrow\RR\cup\{+\infty\}$ is defined by
\begin{equation*}\label{assofunc}
\omega_M(t):=\sup_{p\in\NN}\log\left(\frac{t^p}{M_p}\right)\;\;\;\text{for}\;t>0,\hspace{30pt}\omega_M(0):=0.
\end{equation*}
For an abstract introduction of the associated function we refer to \cite[Chapitre I]{mandelbrojtbook}, see also \cite[Definition 3.1]{Komatsu73}. If $\liminf_{p\rightarrow+\infty}(M_p)^{1/p}>0$, then $\omega_M(t)=0$ for sufficiently small $t$, since $\log\left(\frac{t^p}{M_p}\right)<0\Leftrightarrow t<(M_p)^{1/p}$ holds for all $p\in\NN_{>0}$. Moreover under this assumption $t\mapsto\omega_M(t)$ is a continuous nondecreasing function, which is convex in the variable $\log(t)$ and tends faster to infinity than any $\log(t^p)$, $p\ge 1$, as $t\rightarrow+\infty$. $\lim_{p\rightarrow+\infty}(M_p)^{1/p}=+\infty$ implies that $\omega_M(t)<+\infty$ for any finite $t$ and which shall be considered as a basic assumption for defining $\omega_M$.

Consequently for any normalized $M\in\RR^{\NN}_{>0}$ with $\lim_{p\rightarrow+\infty}(M_p)^{1/p}=+\infty$ its associated weight $\omega_M$ is a normalized weight function and hence satisfying the assumptions from Proposition \ref{Prop13MT88}. In this case one is able to define the so-called log-convex minorant $M^{\on{lc}}=(M^{\on{lc}}_p)_p$ by
\begin{equation}\label{lcminorant}
M^{\on{lc}}_p=\sup_{t>0}\frac{t^p}{\exp(\omega_M(t))},\;\;\;\forall\;p\in\NN,
\end{equation}
see again \cite[Chapitre I]{mandelbrojtbook} and \cite[(3.2)]{Komatsu73}. So $M^{\on{lc}}$ is the largest log-convex sequence satisfying $M^{\on{lc}}\le M$ and we have under theses assumptions on $M$ that $\omega_M=\omega_{M^{\on{lc}}}$.

By definition, for any $t,s>0$ we get
\begin{equation}\label{omegaMspower}
(\omega_M)^s(t)=\omega_M(t^s)=\sup_{p\in\NN}\log\left(\frac{t^{sp}}{M_p}\right)=\sup_{p\in\NN}\log\left(\left(\frac{t^p}{(M_p)^{1/s}}\right)^s\right)=s\omega_{M^{1/s}}(t),
\end{equation}
where $M^{1/s}:=((M_p)^{1/s})_{p\in\NN}$.

We summarize some more well-known facts for this function:

\begin{lemma}\label{assofuncproper}
Let $M\in\hyperlink{LCset}{\mathcal{LC}}$.

\begin{itemize}
\item[$(i)$] $\omega_M$ belongs always to \hyperlink{omset0}{$\mathcal{W}_0$}.

\item[$(ii)$] $\liminf_{p\rightarrow+\infty}(m_p)^{1/p}>0$ implies \hyperlink{om2}{$(\omega_2)$}, $\lim_{p\rightarrow+\infty}(m_p)^{1/p}=+\infty$ implies \hyperlink{om5}{$(\omega_5)$} for $\omega_M$.

\item[$(iii)$] $M$ has \hyperlink{mg}{$(\on{mg})$} if and only if $\omega_M$ has \hyperlink{om6}{$(\omega_6)$}.

\item[$(iv)$] $\omega_M$ has \hyperlink{omnq}{$(\omega_{\on{nq}})$} if and only if $M$ has \hyperlink{mnq}{$(\on{nq})$} and in this case $\lim_{p\rightarrow\infty}(m_p)^{1/p}=+\infty$.

\item[$(v)$] If $M$ satisfies $\hyperlink{beta1}{(\beta_1)}\Leftrightarrow\hyperlink{gamma1}{(\gamma_1)}$, then $\omega_M$ has \hyperlink{omsnq}{$(\omega_{\on{snq}})$} (and which implies \hyperlink{om1}{$(\omega_1)$}).


\item[$(vi)$] If $\omega\in\hyperlink{omset0}{\mathcal{W}_0}$, then $\omega\hyperlink{sim}{\sim}\omega_{W^l}$ for each $l>0$, more precisely we get
    \begin{equation*}\label{goodequivalence}
    \forall\;l>0\;\exists\;C_l>0\;\forall\;t\ge 0:\;\;\;l\omega_{W^l}(t)\le\omega(t)\le 2l\omega_{W^l}(t)+C_l.
    \end{equation*}
\end{itemize}
\end{lemma}

{\itshape Some consequences:} For any $M\in\hyperlink{LCset}{\mathcal{LC}}$ satisfying \hyperlink{gamma1}{$(\gamma_1)$} the function $\omega_M$ is a strong weight, in particular this holds true for any $M\in\hyperlink{SRset}{\mathcal{SR}}$. But in \cite[Section 5]{firstindexpaper} we present a (counter)-example of a strongly log-convex sequence such that $\omega_M$ is a strong weight but \hyperlink{gamma1}{$(\gamma_1)$} is violated. Note that by \cite[Remark 3.3]{sectorialextensions} if $\omega_M$ has \hyperlink{om1}{$(\omega_1)$}, for each $N$ which is equivalent to $M$, we also get that the associated weight functions $\omega_M$ and $\omega_N$ are equivalent. Thus $\omega_N$ is a strong weight too whereas by the characterization in \cite{petzsche} it follows that $N$ does not have \hyperlink{gamma1}{$(\gamma_1)$} either.

\demo{Proof}
$(i)$ See \cite[Definition 3.1]{Komatsu73}.

$(ii)$ That $\lim_{p\rightarrow+\infty}(m_p)^{1/p}=+\infty$ implies \hyperlink{om5}{$(\omega_5)$} for $\omega_M$ follows analogously as $\liminf_{p\rightarrow+\infty}(m_p)^{1/p}>0$ implies \hyperlink{om2}{$(\omega_2)$}, shown in \cite[Lemma 12 $(iv)\Rightarrow(v)$]{BonetMeiseMelikhov07}. Note that by Stirling's formula $\liminf_{p\rightarrow+\infty}(m_p)^{1/p}>0$ is precisely $(M0)$ in \cite{BonetMeiseMelikhov07}.

$(iii)$ See \cite[Proposition 3.6]{Komatsu73}.

$(iv)$ See \cite[Lemma 4.1]{Komatsu73} and \cite[Section 4]{testfunctioncharacterization}.

$(v)$ Follows from \cite[Proposition 4.4]{Komatsu73}.


$(vi)$ See \cite[Lemma 2.5]{sectorialextensions}.
\qed\enddemo

If $M\in\hyperlink{LCset}{\mathcal{LC}}$ and $\Omega=\{W^x: x>0\}$ is denoting the matrix associated to the weight $\omega_M$, then for all $p\in\NN$ we have
\begin{equation}\label{omegaMassofctmatrix}
M_p=\sup_{t\ge0}\frac{t^p}{\exp(\omega_M(t))}=\exp\left(\sup_{t\ge0}\left(p\log(t)-\omega_M(t)\right)\right)=\exp\left(\varphi^{*}_{\omega_M}(p)\right)=W^1_p.
\end{equation}

Let $M\in\RR_{>0}^{\NN}$ (with $M_0=1$) and put
\begin{equation}\label{functionhequ1}
h_M(t):=\inf_{k\in\NN}M_k t^k.
\end{equation}
The functions $h_M$ and $\omega_M$ are related by
\begin{equation*}\label{functionhequ2}
h_M(t)=\exp(-\omega_M(1/t))\;\;\;\forall\;t>0,
\end{equation*}
since $\log(h_M(t))=\inf_{k\in\NN}\log(t^kM_k)=-\sup_{k\in\NN}-\log(t^kM_k)=-\omega_M(1/t)$ (e.g. see also \cite[p. 11]{ChaumatChollet94}). By definition we immediately get:

\begin{lemma}\label{functionhproperties}
Let $M,N\in\RR_{>0}^{\NN}$ be given, then
\begin{itemize}
\item[$(i)$] $t\mapsto h_M(t)$ is nondecreasing,
\item[$(ii)$] $h_{M}(t)\le 1$ for all $t>0$, $h_M(t)=1$ for all $t$ sufficiently large and $\lim_{t\rightarrow 0}h_M(t)=0$,
\item[$(iii)$] $M\le N$ implies $h_M\le h_N$, more generally $M\hyperlink{precsim}{\precsim}N$ implies that $h_M(t)\le h_N(Ct)$ holds for some $C\ge 1$ and all $t>0$,
\item[$(iv)$] for any $s,t>0$ we obtain
\begin{equation*}\label{functionpower}
h_{M^s}(t^s)=(h_M(t))^s.
\end{equation*}
\end{itemize}
\end{lemma}

\subsection{Legendre conjugates of a weight $\omega$}\label{conugateofaweight}
In this section we recall several useful results concerning the conjugates from \cite[Section 3.1]{sectorialextensions}. As we have mentioned there in the study of ultradifferentiable or ultraholomorphic classes defined by weight sequences, the operations of multiplying or dividing the sequence by the factorials play a prominent role. It turns out that when weight functions are considered instead, the corresponding analogous operations are expressed by means of Legendre conjugates, which will be described now. Moreover, the results mentioned in Section \ref{growthindexgamma} below (see also \cite[Section 4.2]{sectorialextensions}, full proofs are given in \cite{firstindexpaper}) will describe how a suitable growth index defined in terms of $\omega$ is increased or decreased by one as it occurs for the analogous growth index introduced in terms of weight sequences by multiplying or dividing by factorials.\vspace{6pt}

In order to prove the extension theorem we will have to deal with $w^l$, where $w^l_j:=\frac{W^l_j}{j!}$. A priori it is not clear that $w^l\in\hyperlink{LCset}{\mathcal{LC}}$ will be valid for some/each $l>0$. But for all $M\in\hyperlink{LCset}{\mathcal{LC}}$ with $\lim_{p\rightarrow+\infty}(m_p)^{1/p}=+\infty$ there exists a close connection between $\omega_m$ and the Legendre conjugate of $\omega_M$, as considered in \cite[Definition 1.4]{PetzscheVogt} and \cite{BonetBraunMeiseTaylorWhitneyextension}, see Lemma \ref{omegaconjugateequivalent} below for more details. This conjugate must not be mixed with $\varphi^{*}_{\omega}$ as considered in \eqref{legendreconjugate}. But since we will always work with the sequences $W^x$ and $w^x$ and their properties introduced in Section \ref{weightmatrixfromfunction} directly we are not going to use $\varphi^{*}_{\omega}$ anymore explicitly.\vspace{6pt}

For any weight function $\omega$ we define its {\itshape upper Legendre conjugate (or envelope)}
\begin{equation*}\label{omegaconjugate}
\omega^{\star}(s):=\sup_{t\ge 0}\{\omega(t)-st\},\;\;\;s\ge 0.
\end{equation*}

We summarize some basic properties listed in \cite[Remark 1.5]{PetzscheVogt} and in \cite[Section 3.1]{sectorialextensions}, for more details we refer also to \cite{firstindexpaper}.

\begin{itemize}
\item[$(i)$] By definition $\omega^{\star}(0)=+\infty$ and, if $\omega$ has in addition \hyperlink{om5}{$(\omega_5)$}, then $\omega^{\star}(s)<+\infty$ for all $s>0$: We have that for any $s>0$ (small) there exists some $C_s>0$ (large) such that for all $t\ge0$ we get $\omega(t)\le st+C_s\Leftrightarrow\omega(t)-st\le C_s$. In this case $\omega^{\star}:(0,+\infty)\rightarrow[0,+\infty)$, $s\mapsto\omega^{\star}(s)$, is nonincreasing, continuous and convex, $\lim_{t\rightarrow+\infty}\omega^{\star}(t)=0$ and $\lim_{t\rightarrow 0}\omega^{\star}(t)=+\infty$.

\item[$(ii)$] Hence for any weight function $\omega$ satisfying \hyperlink{om5}{$(\omega_5)$} we have that $(\omega^{\star})^{\iota}: [0,+\infty)\longrightarrow[0,+\infty)$ is again a (in general not normalized) weight function if we set $(\omega^{\star})^{\iota}(0):=0$ (in particular all assumptions in Proposition \ref{Prop13MT88} are satisfied).
\end{itemize}

For any $h:(0,+\infty)\rightarrow[0,+\infty)$ which is nonincreasing and such that $\lim_{s\rightarrow 0}h(s)=+\infty$, we can define the so-called {\itshape lower Legendre conjugate (or envelope)} $h_{\star}:[0,+\infty)\rightarrow[0,+\infty)$ of $h$ by
\begin{equation}\label{omegaconjugate0}
h_{\star}(t):=\inf_{s>0}\{h(s)+ts\},\hspace{15pt}t\ge 0.
\end{equation}
We are summarizing some facts for this conjugate, see also \cite[Section 3.1]{sectorialextensions}.

$h_{\star}$ is clearly nondecreasing, continuous and concave, and $\lim_{t\rightarrow+\infty}h_{\star}(t)=+\infty$, see \cite[(8), p. 156]{Beurling72}. Moreover, if $\lim_{s\to +\infty}h(s)=0$ then $h_{\star}(0)=0$, and so $h_{\star}$ is a weight function.

In this paper this conjugate will be mainly applied to the case $h(t):=\omega^{\iota}(t)=\omega(1/t)$, where $\omega$ is a weight function, so that $(\omega^{\iota})_{\star}$ is again a weight function; in particular, we will frequently find the case $h(t)=\omega^{\iota}_M(t)=\omega_M(1/t)$ for $M\in\RR_{>0}^{\NN}$ with $\lim_{p\rightarrow+\infty}(M_p)^{1/p}=+\infty$. In case $\omega$ is a weight function satisfying \hyperlink{om5}{$(\omega_5)$}, $\sigma:=(\omega^{\star})_{\star}$ is a weight function and it is indeed the least concave majorant of $\omega$ in the sense that, if $\tau:[0,+\infty)\rightarrow[0,+\infty)$ is concave and $\omega\le\tau$, then $\sigma=(\omega^{\star})_{\star}\le\tau$, see ~\cite{roever}.\vspace{6pt}

The proof of \cite[Lemma 5.1]{sectorialextensions} actually shows that $\alpha\hyperlink{sim}{\sim}\beta$ implies $(\alpha^{\iota})_{\star}\hyperlink{sim}{\sim}(\beta^{\iota})_{\star}$.\vspace{6pt}

We recall now the following two useful results which establish a relation between $\omega^{\star}$ (respectively $\omega_{W^x}^{\star}$) and $\omega_{w^x}$ (and hence $h_{w^x}$), for the proofs we refer to \cite[Lemma 3.1, (3.5), Corollary 3.2]{sectorialextensions} (showing in \cite[Lemma 3.1]{sectorialextensions} even a more general version of \eqref{Dynkinequiv1} by involving a parameter) respectively see also \cite[Lemma 3.10, Corollary 3.11]{whitneyextensionmixedweightfunction}.

\begin{lemma}\label{omegaconjugateequivalent}
Let $M\in\RR_{>0}^{\NN}$ such that $\lim_{p\rightarrow+\infty}(m_p)^{1/p}=+\infty$, then
\begin{equation}\label{Dynkinequiv1}
\forall\;s>0:\;\;\;\omega_{M}^{\star}(s)\le\omega_{m}\left(\frac{1}{s}\right)\le\omega_{M}^{\star}\left(\frac{s}{e}\right)\Leftrightarrow(\omega_M^{\star})^{\iota}(s)\le\omega_m(s)\le(\omega_M^{\star})^{\iota}(se).
\end{equation}
\end{lemma}

Combining this result with the approach from \cite[Section 5]{compositionpaper} (see also \cite[Lemma 2.5]{sectorialextensions}) we get the following consequences for the sequences $W^x\in\Omega$.

\begin{corollary}\label{omegaconjugateequivalentcor}
Let $\omega\in\hyperlink{omset0}{\mathcal{W}_0}$ be given with \hyperlink{om5}{$(\omega_5)$} and $\Omega=\{W^x: x>0\}$ its associated weight matrix, then
\begin{equation*}\label{omegaconjugateequivalentcorequ1}
\forall\;x>0\;\exists\;C_x\ge 1\;\forall\;s>0:\;\;\;x\omega_{W^x}^{\star}(\frac{s}{x})\le\omega^{\star}(s)\le 2x\omega_{W^x}^{\star}(\frac{s}{2x})+C_x
\end{equation*}
and
\begin{equation*}\label{omegaconjugateequivalentcorequ2}
\forall\;x>0\;\exists\;C_x\ge 1\;\forall\;s>0:\;\;\;x\omega_{w^x}\left(\frac{x}{es}\right)\le\omega^{\star}(s)\le 2x\omega_{w^x}\left(\frac{2x}{s}\right)+C_x,
\end{equation*}
respectively equivalently
\begin{equation}\label{omegaconjugateequivalentcorequ3}
\forall\;x>0\;\exists\;C_x\ge 1\;\forall\;s>0:\;\;\;h_{w^x}\left(\frac{es}{x}\right)^x\ge\exp(-\omega^{\star}(s))\ge \exp(-C_x)h_{w^x}\left(\frac{s}{2x}\right)^{2x}.
\end{equation}
\end{corollary}

Generalizing \cite[Proposition 3.6]{Komatsu73} we have shown in \cite[Proposition 3.6]{testfunctioncharacterization} respectively in \cite[Remark 2.5]{whitneyextensionweightmatrix} the following characterization.

\begin{lemma}\label{functionh1}
Let $M,N\in\hyperlink{LCset}{\mathcal{LC}}$ be given. Then $\hyperlink{genmg}{(M,N)_{(\on{mg})}}$ holds if and only if

\begin{equation}\label{equfunctionh1}
\exists\;A\ge 1\;\forall\;t>0:\;\;h_{M}(t)\le h_{N}(At)^2\Leftrightarrow 2\omega_N(t)\le\omega_M(At).
\end{equation}
\end{lemma}

Note that for the proof that $\hyperlink{genmg}{(M,N)_{(\on{mg})}}$ implies \eqref{equfunctionh1} it is sufficient to assume $M\in\RR_{>0}^{\NN}$ and $N\in\hyperlink{LCset}{\mathcal{LC}}$, see \cite[Lemma 3.5]{Komatsu73} and the proof of \cite[Proposition 3.6]{testfunctioncharacterization}. By using the conjugate and \eqref{Dynkinequiv1} we can translate this mixed condition to $\omega^{\star}$ respectively to $h_{w^l}$ as well, see also \cite[Lemma 3.13]{whitneyextensionmixedweightfunction} where we have started with two sequences.

\begin{lemma}\label{functionh2}
Let $\omega\in\hyperlink{omset0}{\mathcal{W}_0}$ with \hyperlink{om5}{$(\omega_5)$} be given and $\Omega=\{W^l: l>0\}$ its associated weight matrix. Then
\begin{equation}\label{functionh2equ0}
\forall\;l>0\;\forall\;s>0:\;\;\;2\omega^{\star}_{W^{2l}}(s)\le\omega^{\star}_{W^l}(2s),
\end{equation}
and
\begin{equation}\label{functionh2equ1}
\exists\;A\ge 1\;\forall\;l>0\;\forall\;s>0:\;\;\;h_{w^l}(s)\le h_{w^{2l}}(As)^2.
\end{equation}
\end{lemma}

\demo{Proof}
By \eqref{newmoderategrowth} in \ref{weightmatrixfromfunction} and \eqref{equfunctionh1} for all $l>0$ and all $t>0$ we get $2\omega_{W^{2l}}(t)\le\omega_{W^l}(t)$ and both mappings $s\mapsto\omega^{\star}_{W^{2l}}(s)$ and $s\mapsto\omega^{\star}_{W^l}(s)$ are well-defined, see also $(ii)$ in Lemma \ref{omegaconjugateequivalent}. Hence
$$2\omega^{\star}_{W^{2l}}(s)=\sup_{t\ge 0}\{2\omega_{W^{2l}}(t)-(2s)t\}\le\sup_{t\ge 0}\{\omega_{W^l}(t)-(2s)t\}=\omega^{\star}_{W^l}(2s),$$
which proves the first part.

For the proof of \eqref{functionh2equ1} we combine \eqref{Dynkinequiv1} and \eqref{functionh2equ0} to get
$$2\omega_{w^{2l}}(1/s)\le 2\omega_{W^{2l}}^{\star}(s/e)\le\omega^{\star}_{W^l}(2s/e)\le\omega_{w^l}(e/(2s)),$$
hence $h_{w^{2l}}(s)^2=\exp(-2\omega_{w^{2l}}(1/s))\ge\exp(-\omega_{w^l}(e/(2s)))=h_{w^l}(2s/e)$ follows.
\qed\enddemo

\section{Growth index $\gamma(\omega)$}\label{growthindexgamma}
In this section, for a given weight function $\omega$, we recall the definition of $\gamma(\omega)$ from \cite[Section 4.2]{sectorialextensions} and state some properties which will be needed in the main results below. More details, the relation to the growth index $\gamma(M)$ introduced by V. Thilliez in \cite{Thilliezdivision} and precise proofs are given in \cite{firstindexpaper}. At this point we want to give some motivation how to come up with the definition given in \cite[Section 4.2]{sectorialextensions}.\vspace{6pt}

Recall that for any given $M\in\RR_{>0}^{\NN}$ and $t,s>0$ we have $\omega_{M^s}(t)=s\omega_M(t^{1/s})$, see \eqref{omegaMspower}. Moreover clearly $L\hyperlink{approx}{\approx}M\Leftrightarrow L^s\hyperlink{approx}{\approx}M^s$ for any $s>0$.\vspace{6pt}

By Proposition \ref{Prop13MT88} condition \hyperlink{omsnq}{$(\omega_{\on{snq}})$} is equivalent to having $\limsup_{t\rightarrow+\infty}\frac{\omega(Kt)}{\omega(t)}<K$ for some $K>1$. The case $\omega=\omega_M$ yields that for any $s>0$ and $A>1$ by putting $t'=t^{1/s}$ we get
$$\frac{\omega_{M^s}(At)}{\omega_{M^s}(t)}=\frac{\omega_M(A^{1/s}t^{1/s})}{\omega_M(t^{1/s})}=\frac{\omega_M(A^{1/s}t')}{\omega_M(t')},$$
see \eqref{omegaMspower}. This shows that for any $s>1$ the weight $\omega_{M^s}$ has automatically \hyperlink{omsnq}{$(\omega_{\on{snq}})$} if one assumes that $\omega_M$ has this property: take $A:=K^s$ and note that $A>K$ holds.

But for the case $0<s<1$ this is not true in general, however it holds if we assume that $\omega_M$ satisfies
\begin{equation*}\label{newindex0}
\exists\;K>1:\;\;\;\limsup_{t\rightarrow+\infty}\frac{\omega_M(K^{1/s}t)}{\omega_M(t)}<K\Longleftrightarrow\;\exists\;K>1:\;\;\;\limsup_{t\rightarrow\infty}\frac{\omega_M(Kt)}{\omega_M(t)}<K^s(<K).
\end{equation*}
This motivates the following definition: Let $\omega$ be a weight function and $\gamma>0$, then introduce
\begin{equation*}\label{newindex1}
(P_{\omega,\gamma}):\Longleftrightarrow\exists\;K>1:\;\;\;\limsup_{t\rightarrow+\infty}\frac{\omega(K^{\gamma}t)}{\omega(t)}<K.
\end{equation*}
Note: If $(P_{\omega,\gamma})$ holds for some $K>1$, then also $(P_{\omega,\gamma'})$ is satisfied for all $\gamma'\le\gamma$ with the same $K$. Moreover we restrict ourselves to $\gamma>0$, because for $\gamma\le 0$ condition $(P_{\omega,\gamma})$ is satisfied for all weights $\omega$ (since $\omega$ is nondecreasing and $K>1$).

Finally we put
\begin{equation}\label{newindex2}
\gamma(\omega):=\sup\{\gamma>0: (P_{\omega,\gamma})\;\;\text{is satisfied}\}.
\end{equation}
So for any $s<\gamma(\omega)$ the weight $\omega^s:t\mapsto\omega(t^s)$ has property \hyperlink{omsnq}{$(\omega_{\on{snq}})$}. Let $\omega,\sigma$ satisfy $\sigma\hyperlink{sim}{\sim}\omega\Leftrightarrow\sigma^s\hyperlink{sim}{\sim}\omega^s$, $s>0$, then $\gamma(\sigma)=\gamma(\omega)$: For this note that each $(P_{\cdot,\gamma})$ is stable with respect to $\hyperlink{sim}{\sim}$ since \hyperlink{omsnq}{$(\omega_{\on{snq}})$} is clearly stable with respect to this relation. By definition and \eqref{omegaMspower} we immediately get
\begin{equation}\label{newindex3}
\forall\;s>0:\;\;\;\gamma(\omega^{1/s})=s\gamma(\omega),\hspace{30pt}\forall\;s>0\;\forall\;M\in\RR_{>0}^{\NN}:\;\;\;\gamma(\omega_{M^s})=\gamma((\omega_M)^{1/s})=s\gamma(\omega_M).
\end{equation}
Consequently, if $M\in\hyperlink{LCset}{\mathcal{LC}}$, then each $M^{1/s}$ too and for each $s<\gamma(\omega_M)$ the sequence $M^{1/s}$ satisfies \hyperlink{mnq}{$(\on{nq})$}, see \cite[Lemma 4.1]{Komatsu73}.

In \cite[Lemma 4.3]{sectorialextensions} we have shown the following characterization.

\begin{lemma}\label{strongweightspace}
A weight function $\omega$ satisfies \hyperlink{omsnq}{$(\omega_{\on{snq}})$} if and only if $\gamma(\omega)>1$.
\end{lemma}

Summarizing everything from above we have shown that for each $\omega$ we have $\{\gamma\ge 0: (P_{\omega,\gamma})\;\text{holds}\}=[0,\gamma(\omega))$. Note that $(P_{\omega,0})$ is trivially satisfied.\vspace{6pt}

The proofs of all further results in this section are given in full detail in \cite[Section 2]{firstindexpaper} (already stated in \cite[Section 4.2]{sectorialextensions}).

\begin{lemma}\label{omega1weightspace}
A weight function $\omega$ satisfies \hyperlink{om1}{$(\omega_1)$} if and only if $\gamma(\omega)>0$.
\end{lemma}

\begin{lemma}\label{growthindexgammalemma}
Let $\omega$ be given. If $\gamma(\omega)>1$, then the function $(\omega^{\star})^{\iota}$ has \hyperlink{om1}{$(\omega_1)$} and we obtain
$$\gamma(\omega)=\gamma((\omega^{\star})^{\iota})+1.$$
In particular $\gamma(\omega)>2$ implies that $(\omega^{\star})^{\iota}$ satisfies \hyperlink{omsnq}{$(\omega_{\on{snq}})$}, too.
\end{lemma}

Using the previous result we show in \cite{firstindexpaper} the following consequence.

\begin{corollary}\label{growthindexgammacorollary}
Let $M\in\hyperlink{LCset}{\mathcal{LC}}$ be given satisfying $\gamma(\omega_M)>1$. Then we obtain
\begin{equation*}\label{growthindexgammacorollaryequ1}
\gamma(\omega_M)=\gamma(\omega_m)+1,
\end{equation*}
and consequently $\omega_m\in\hyperlink{omset1}{\mathcal{W}}$.

In particular $\gamma(\omega_M)>2$ implies that $\omega_m$ has \hyperlink{omsnq}{$(\omega_{\on{snq}})$}.\vspace{6pt}

If $\omega\in\hyperlink{omset0}{\mathcal{W}_0}$ with $\gamma(\omega)>1$ (i.e. $\omega$ is a strong weight by Lemma \ref{strongweightspace}), then
\begin{equation}\label{growthindexgammacorollaryequ2}
\forall\;x>0:\;\;\;\gamma(\omega)=\gamma(\omega_{W^x})=\gamma(\omega_{w^x})+1=\gamma((\omega^{\star})^{\iota})+1.
\end{equation}
\end{corollary}





On the one hand Lemma \ref{growthindexgammalemma} proves $\gamma(\omega)=\gamma((\omega^{\star})^{\iota})+1$ provided $\gamma(\omega)>1$ and the next two results illustrate the effect of applying the {\itshape lower Legendre envelope} defined in \eqref{omegaconjugate0}.

\begin{lemma}\label{lowerenvelop}
Let $\omega$ be given with $\gamma(\omega)>0$. Then we obtain
$$\gamma(\omega)+1=\gamma((\omega^{\iota})_{\star}).$$
\end{lemma}




Finally, the next result is the converse statement of Lemma \ref{growthindexgammalemma} and Corollary \ref{growthindexgammacorollary}.

\begin{lemma}\label{lowerenvelop1}
Let $M\in\hyperlink{LCset}{\mathcal{LC}}$ with $(m_p)^{1/p}\rightarrow+\infty$ as $p\rightarrow+\infty$ and $\gamma(\omega_m)>0$. We set $L=(L_p)_p$, $L_p:=p!m_p^{\on{lc}}$ (see \eqref{lcminorant}) and obtain
\begin{equation*}\label{lowerenvelopequ1}
\gamma(\omega_m)+1=\gamma(\omega_L).
\end{equation*}
On the other hand, let $\omega\in\hyperlink{omset0}{\mathcal{W}_0}$ be given with $\gamma(\omega)>0$ (equivalently $\omega\in\hyperlink{omset1}{\mathcal{W}}$), $\Omega=\{W^x: x>0\}$ the associated weight matrix and put $\widehat{W}^x:=(p!W^x_p)_{p\in\NN}$. Then we obtain
\begin{equation}\label{lowerenvelopequ2}
\forall\;x>0:\;\;\;\gamma(\omega_{W^x})+1=\gamma(\omega)+1=\gamma((\omega^{\iota})_{\star})=\gamma((\omega_{W^x}^{\iota})_{\star})=\gamma(\omega_{\widehat{W}^x}).
\end{equation}
\end{lemma}




\section{Existence of sectorially flat functions in the weight function setting}\label{sectoriallyflat}
\subsection{Construction of outer functions}\label{constructionouterfunction}
Let $\omega$ be a given weight function satisfying \hyperlink{om5}{$(\omega_5)$}. As we have pointed out in Section \ref{conugateofaweight} it follows that $(\omega^{\star})^{\iota}$ is a weight function too (in the sense of \ref{weightfunctionclasses} and \cite[Section 2.2]{sectorialextensions}). Thus it is possible to apply directly the proofs from \cite[Section 6.1]{sectorialextensions} and replace there $\tau$ by $(\omega^{\star})^{\iota}$ and $\tau^{\iota}$ by $\omega^{\star}$. Hence we obtain the following results which should be compared with the approach in \cite{sectorialextensions} and the results obtained by V. Thilliez in \cite[Section 2.1]{Thilliezdivision}.

The difference to the previous works is that on the one hand we do not start with a weight sequence $M$, consider the associated function $h_m$ and assuming some growth control in terms of $M$ (as in \cite{Thilliezdivision}). And on the other hand we also do not start with a weight function $\omega$ and assuming some growth control on it expressed in terms of $\gamma(\omega)$ directly. In this new approach we start with a weight function $\omega$ and assume that we have some control in terms of $\gamma((\omega^{\star})^{\iota})$ and so the conjugate from Section \ref{conugateofaweight} is involved in our considerations.\vspace{6pt}

First we obtain a characterization of \hyperlink{omsnq}{$(\omega_{\text{snq}})$} for $(\omega^{\star})^{\iota}$, compare this with \cite[Lemma 2.1.1]{Thilliezdivision} and \cite[Lemma 6.1]{sectorialextensions}.

\begin{lemma}\label{functionh3}
Let $\omega$ be a weight function with \hyperlink{om5}{$(\omega_5)$}. Then we have $\gamma((\omega^{\star})^{\iota})>1$ (equivalently \hyperlink{omsnq}{$(\omega_{\on{snq}})$}) if and only if
\begin{equation*}\label{functionh3equ}
\exists\;C\ge 1\;\forall\;y>0:\;\;\;\int_0^1-\omega^{\star}(ty)dt\ge-C(\omega^{\star}(y)+1).
\end{equation*}
\end{lemma}

The next result is analogous to \cite[Lemma 2.1.2]{Thilliezdivision} respectively \cite[Lemma 6.2]{sectorialextensions}.

\begin{lemma}\label{Lemma212}
Let $\omega$ be a weight function with \hyperlink{om5}{$(\omega_5)$} and such that $\gamma((\omega^{\star})^{\iota})>1$, then
$$\int_{-\infty}^{+\infty}\frac{-\omega^{\star}(|t|)}{1+t^2}dt>-\infty.$$
\end{lemma}

Finally we transfer \cite[Lemma 2.1.3]{Thilliezdivision} to the weight function case, see \cite[Lemma 6.3]{sectorialextensions}.

\begin{lemma}\label{Lemma213}
Let $\omega$ be a weight function with \hyperlink{om5}{$(\omega_5)$} and such that $\gamma((\omega^{\star})^{\iota})>1$. Then for all $a>0$ the function $F_a$ defined by
$$F_a(w):=\exp\left(\frac{1}{\pi}\int_{-\infty}^{+\infty} \frac{-a\omega^{\star}(|t|)}{1+t^2}\frac{itw-1}{it-w}dt\right),$$
is holomorphic in the right half-plane $H_1:=\{\xi\in\CC: \Re(\xi)>0\}$, and there exist constants $A,B\ge 1$ (large) depending only on $\omega$ such that
\begin{equation}\label{Lemma213equ}
\forall\;\xi\in H_1:\;\;\;B^{-a}\exp(-2a\omega^{\star}(B^{-1}\Re(\xi)))\le|F_a(\xi)|\le \exp\left(-\frac{a}{2}\omega^{\star}(A|\xi|)\right).
\end{equation}
\end{lemma}

\subsection{Construction of optimal sectorially flat functions}\label{sectoriallyflatfunction}
Using the results from Section \ref{constructionouterfunction} the aim is now to transfer \cite[Theorem 2.3.1]{Thilliezdivision} to the weight function setting. For this we can again partially follow the proof of \cite[Theorem 6.7]{sectorialextensions} by replacing $\tau$ by $(\omega^{\star})^{\iota}$ and $\tau^{\iota}$ by $\omega^{\star}$.

To state the following theorem it is necessary to have $\gamma((\omega^{\star})^{\iota})>0$. Note that by combining Lemma \ref{strongweightspace} and Lemma \ref{growthindexgammalemma} we see that if a weight function $\omega$ satisfies \hyperlink{omsnq}{$(\omega_{\on{snq}})$}, then $\gamma((\omega^{\star})^{\iota})=\gamma(\omega)-1>0$ follows. On the other hand, as mentioned in Lemma \ref{omega1weightspace} for any weight function $\omega$ we have $\gamma(\omega)>0$ if and only if \hyperlink{om1}{$(\omega_1)$} holds true.\vspace{6pt}

From now on we will frequently assume that $\omega\in\hyperlink{omset1}{\mathcal{W}}$ since in the proofs we wish to make use of the representation \eqref{equaEqualitySpacesWeightFunctionMatrix} and work with the sequences in the associated weight matrix $\Omega=\{W^x: x>0\}$ by applying \eqref{newexpabsorb}.

\begin{theorem}\label{Theorem231}
Let $\omega\in\hyperlink{omset1}{\mathcal{W}}$ be given with \hyperlink{om5}{$(\omega_5)$} and such that $\gamma((\omega^{\star})^{\iota})>0$. Then for any $0<\gamma<\gamma((\omega^{\star})^{\iota})$ there exist constants $K_1,K_2,K_3>0$ depending only on $\omega$ and $\gamma$ such that for all $a>0$ there exists a function $G_a\in\mathcal{A}_{\{W^{y(a)}\}}(S_{\gamma})\subseteq\mathcal{A}_{\{\omega\}}(S_{\gamma})$ satisfying
\begin{equation}\label{Theorem231equ}
\forall\;\xi\in S_{\gamma}:\;\;K_1^{-a}\exp(-2a\omega^{\star}(K_2|\xi|))\le|G_a(\xi)|\le\exp\left(-\frac{a}{2}\omega^{\star}(K_3|\xi|)\right).
\end{equation}
\end{theorem}

\demo{Proof}
For the proof of \eqref{Theorem231equ} we apply \cite[Theorem 6.7]{sectorialextensions} with $(\omega^{\star})^{\iota}$ instead of $\tau$, $\omega^{\star}$ instead of $\tau^{\iota}$.\vspace{6pt}

It remains to show that $G_a\in\mathcal{A}_{\{W^y\}}(S_{\gamma})$ for some $y=y(a)>0$, the inclusion $\mathcal{A}_{\{W^{y(a)}\}}(S_{\gamma})\subseteq\mathcal{A}_{\{\omega\}}(S_{\gamma})$ holds true by \hyperlink{om1}{$(\omega_1)$}. We follow \cite[Theorem 2.3.1]{Thilliezdivision} and recall briefly the proof of the right-hand side of \eqref{Theorem231equ}.

Let $a>0$ be arbitrary, take $s,\delta>0$ such that $\gamma<\delta<\gamma((\omega^{\star})^{\iota})$, $s\delta<1<s\gamma((\omega^{\star})^{\iota})$ and so $\gamma(((\omega^{\star})^{\iota})^{1/s})>1$ (see \eqref{newindex3}). We apply Lemma \ref{Lemma213} to the weight $((\omega^{\star})^{\iota})^{1/s}$ instead of $(\omega^{\star})^{\iota}$ and note that $(((\omega^{\star})^{\iota})^{1/s})^{\iota}(t)=((\omega^{\star})^{\iota})^{1/s}(1/t)=(\omega^{\star})^{\iota}(1/t^{1/s})=\omega^{\star}(t^{1/s})=(\omega^{\star})^{1/s}(t)$.
Then put
$$G_a(\xi)=F_a(\xi^s),\;\;\;\;\xi\in S_{\delta},$$
where $F_a$ is the function from Lemma \ref{Lemma213}.

Let $A$ be the constant coming from the right-hand side of \eqref{Lemma213equ} applied to this situation, so
\begin{align*}
|G_a(\xi)|&=|F_a(\xi^s)|\le \exp\left(-\frac{a}{2}(\omega^{\star})^{1/s}(A|\xi|^s)\right)=\exp\left(-\frac{a}{2}(\omega^{\star})^{1/s}((A^{1/s}|\xi|)^s)\right)
\\&
=\exp\left(-\frac{a}{2}\omega^{\star}(A^{1/s}|\xi|)\right).
\end{align*}
First put in the estimate above $A_1:=A^{1/s}$ and then, by using \eqref{omegaconjugateequivalentcorequ3} for any $y>0$, we get
$$\exp\left(-\frac{a}{2}\omega^{\star}(A_1|\xi|)\right)\le h_{w^y}\left(\frac{eA_1|\xi|}{y}\right)^{ya/2}.$$
Hence taking $y:=2a^{-1}$ proves $|G_a(\xi)|\le h_{w^y}\left(\frac{aeA_1|\xi|}{2}\right)$.

Then take $\varepsilon>0$ such that $\varepsilon<\min\{1,\delta-\gamma\}\pi/2$, so for any $\xi\in S_{\gamma}$ the closed disc centered at $\xi$ and radius $\sin(\varepsilon)|\xi|$ belongs to $S_{\delta}$. The Cauchy formula applied to this disc denoted by $B_{\varepsilon}$ and the right hand side of \eqref{Theorem231equ} applied to $S_{\delta}$ imply
\begin{equation}\label{Theorem231equ1}
\forall\;\xi\in B_{\varepsilon}\;\forall\;j\in\NN:\;\;\;|G_a^{(j)}(\xi)|\le \frac{j!}{(\sin(\varepsilon)|\xi|)^j}h_{w^y}(A_2(1+\sin(\varepsilon))|\xi|).
\end{equation}
By definition \eqref{functionhequ1} it is clear that $h_{w^y}(A_2(1+\sin(\varepsilon))|\xi|)\le(A_2(1+\sin(\varepsilon))|\xi|)^jw^y_j$ holds for all $j\in\NN$, hence for all $j\in\NN$ and $\xi\in S_{\gamma}$:
\begin{equation*}\label{Theorem231equ2}
|G_a^{(j)}(\xi)|\le\frac{j!}{(\sin(\varepsilon)|\xi|)^j}(A_2(1+\sin(\varepsilon))|\xi|)^jw^y_j=\left(\frac{A_2(1+\sin(\varepsilon))}{\sin(\varepsilon)}\right)^jW^y_j.
\end{equation*}
\qed\enddemo

Now we can transfer \cite[Lemma 2.3.2]{Thilliezdivision} to the weight function case.

\begin{lemma}\label{Lemma232}
Let $\omega\in\hyperlink{omset1}{\mathcal{W}}$ be given with \hyperlink{om5}{$(\omega_5)$} and such that $\gamma((\omega^{\star})^{\iota})>0$, $\Omega=\{W^x: x>0\}$ the matrix associated with $\omega$. Then the derivatives of $G_a\in\mathcal{A}_{\{\omega\}}(S_{\gamma})$ and $\frac{1}{G_a}$ constructed in Theorem \ref{Theorem231} satisfy the following estimates:
\begin{equation}\label{Lemma232equ}
\forall\;a>0\;\exists\;E_1,E_2>0\;\forall\;\xi\in S_{\gamma}\;\forall\;j\in\NN:\;\;\;|G_a^{(j)}(\xi)|\le E_1^{j+1}W^{4/a}_jh_{w^{4/a}}(E_2|\xi|),
\end{equation}
\begin{equation}\label{Lemma232equ1}
\forall\;y>0\;\exists\;x>0\;\exists\;a>0\;\exists\;E_3,E_4,E_5>0\;\forall\;\xi\in S_{\gamma}\;\forall\;j\in\NN:\;\;\;\left|\left(\frac{1}{G_a}\right)^{(j)}(\xi)\right|\le E_3 E_4^jW^{x}_j\frac{1}{h_{w^y}(E_5|\xi|)}.
\end{equation}
\end{lemma}

By $(i)$ in Lemma \ref{functionhproperties} we see that in \eqref{Lemma232equ} we can choose $E_1=E_2$ and in \eqref{Lemma232equ1} we can choose $E_3=E_4=E_5^{-1}$.

\demo{Proof}
To prove \eqref{Lemma232equ} we start in the estimate with \eqref{Theorem231equ1} and then use \eqref{functionh2equ1} in Lemma \ref{functionh2} with the universal constant $A$ there and finally the definition of \eqref{functionhequ1} to obtain for any $\xi\in S_{\gamma}$:
\begin{align*}
|G_a^{(j)}(\xi)|&\le\frac{j!}{(\sin(\varepsilon)|\xi|)^j}h_{w^{2/a}}(A_2(1+\sin(\varepsilon))|\xi|)\le \frac{j!}{(\sin(\varepsilon)|\xi|)^j}(h_{w^{4/a}}(A_2A(1+\sin(\varepsilon))|\xi|))^2
\\&
\le\frac{j!}{(\sin(\varepsilon)|\xi|)^j}(A_2A(1+\sin(\varepsilon))|\xi|)^jw_j^{4/a}h_{w^{4/a}}(A_2A(1+\sin(\varepsilon))|\xi|)
\\&
=\left(\frac{A_2A(1+\sin(\varepsilon))}{\sin(\varepsilon)}\right)^jW^{4/a}_jh_{w^{4/a}}(A_2A(1+\sin(\varepsilon))|\xi|).
\end{align*}
For \eqref{Lemma232equ1}, first we use the left hand side of \eqref{Theorem231equ} on the sector $S_{\delta}$ and the Cauchy formula on the closed disc $B_{\varepsilon}$ from the proof of Theorem \ref{Theorem231}. Let $a,b,x,y>0$ be arbitrary, then for all $j\in\NN$ and $\xi\in S_{\gamma}$ we get by applying \eqref{omegaconjugateequivalentcorequ3} twice
\begin{align*}
\left|\left(\frac{1}{G_a}\right)^{(j)}(\xi)\right|&\le\frac{j!}{K_1^{-a}(\sin(\varepsilon)|\xi|)^j}\exp(2a\omega^{\star}(K_2(1-\sin(\varepsilon))|\xi|))
\\&
=K_1^a\frac{j!}{(\sin(\varepsilon)|\xi|)^j}\frac{\exp(2a\omega^{\star}(K_2(1-\sin(\varepsilon))|\xi|))}{\exp(-b\omega^{\star}(K_2(1-\sin(\varepsilon))|\xi|))}\exp(-b\omega^{\star}(K_2(1-\sin(\varepsilon))|\xi|))
\\&
=K_1^a\frac{j!}{(\sin(\varepsilon)|\xi|)^j}\exp((2a+b)\omega^{\star}(K_2(1-\sin(\varepsilon))|\xi|))\exp(-b\omega^{\star}(K_2(1-\sin(\varepsilon))|\xi|))
\\&
\le K_1^a\frac{j!}{(\sin(\varepsilon)|\xi|)^j}\exp((2a+b)\omega^{\star}(K_2(1-\sin(\varepsilon))|\xi|))h_{w^x}\left(\frac{eK_2(1-\sin(\varepsilon))|\xi|}{x}\right)^{xb}
\\&
\le K_1^a\frac{j!}{(\sin(\varepsilon)|\xi|)^j}\frac{1}{C^{2a+b}_yh_{w^y}\left(\frac{K_2(1-\sin(\varepsilon))|\xi|}{2y}\right)^{2y(2a+b)}}h_{w^x}\left(\frac{eK_2(1-\sin(\varepsilon))|\xi|}{x}\right)^{xb}.
\end{align*}
Let $y>0$ be arbitrary but from now on fixed and we choose $a:=\frac{1}{8y}$, $b:=\frac{1}{4y}$ and finally $x:=b^{-1}=4y$. By definition of $h_{w^x}$ in \eqref{functionhequ1} we have for all $j\in\NN$ that with these choices $h_{w^x}\left(\frac{eK_2(1-\sin(\varepsilon))|\xi|)}{x}\right)^{xb}=h_{w^x}\left(\frac{eK_2(1-\sin(\varepsilon))|\xi|)}{b^{-1}}\right)\le (beK_2(1-\sin(\varepsilon))|\xi|)^jw^x_j$. Moreover we get $2y(2a+b)=1$, $2a+b=(2y)^{-1}$, hence
\begin{align*}
\left|\left(\frac{1}{G_a}\right)^{(j)}(\xi)\right|&\le K_1^a\frac{j!}{(\sin(\varepsilon)|\xi|)^j}\frac{1}{C^{2a+b}_yh_{w^y}((2a+b)K_2(1-\sin(\varepsilon))|\xi|)}(beK_2(1-\sin(\varepsilon))|\xi|)^jw^x_j
\\&
=\frac{K_1^a}{(C_y)^{1/(2y)}}\left(\frac{beK_2(1-\sin(\varepsilon))}{\sin(\varepsilon)}\right)^jW^x_j\frac{1}{h_{w^y}(1/(2y)K_2(1-\sin(\varepsilon))|\xi|)}.
\end{align*}
\qed\enddemo



\section{The Roumieu type ultraholomorphic extension theorem in the weight function setting}\label{ultraholomorphicextension}
Using the functions $G_a$ constructed in Section \ref{sectoriallyflatfunction} we are now able to transfer \cite[Theorem 3.2.1]{Thilliezdivision} to classes of ultraholomorphic functions defined by strong weight functions $\omega$. Since in the original proof for strongly regular weight sequences several tools from ultradifferenetiable functions have been used, we also recall briefly the definitions for such classes now.

\subsection{Ultradifferentiable classes $\mathcal{E}_{\{\omega\}}$ and $\mathcal{E}_{\{M\}}$}\label{ultradifferentiableclasses}
Let $d\in\NN_{>0}$, $U\subseteq\RR^d$ be nonempty open, $M\in\RR_{>0}^{\NN}$ and $\omega\in\hyperlink{omset0}{\mathcal{W}_0}$. The space of ultradifferentiable functions of (global) Roumieu type in terms of a weight sequence $M$ is defined by
\begin{equation}\label{defofm}
\mathcal{E}_{\{M\}}(U,\CC):=\{f\in\mathcal{E}(U,\CC):\;\;\exists\;h>0:\;\|f\|_{M,U,h}<+\infty\},
\end{equation}
and in terms of a weight function $\omega$ by
\begin{equation}\label{defofomega}
\mathcal{E}_{\{\omega\}}(U,\CC):=\{f\in\mathcal{E}(U,\CC):\;\;\exists\;l>0:\;\|f\|_{\omega,U,l}<+\infty\}
\end{equation}
where
\begin{equation*}\label{semi-norm-1}
\|f\|_{M,U,h}:=\sup_{k\in\NN^d,x\in U}\frac{|f^{(k)}(x)|}{h^{|k|} M_{|k|}}
\end{equation*}
respectively
\begin{equation*}\label{semi-norm-2}
\|f\|_{\omega,U,l}:=\sup_{k\in\NN^d,x\in U}\frac{|f^{(k)}(x)|}{\exp(\frac{1}{l}\varphi^{*}_{\omega}(l|k|))}.
\end{equation*}
Here $f^{(k)}$ is standing for taking partial derivatives of $f$ with respect to $k=(k_1,\dots k_d)\in\NN^d$. The classes are endowed with the locally convex topology given by
\begin{equation}\label{repr1}
\mathcal{E}_{\{M\}}(U,\CC):=\underset{h>0}{\varinjlim}\;\mathcal{E}_{M,h}(U,\CC),\hspace{25pt}\mathcal{E}_{\{\omega\}}(U,\CC):=\underset{l>0}{\varinjlim}\;\mathcal{E}_{\omega,l}(U,\CC),
\end{equation}
where $\mathcal{E}_{M,h}(U,\CC):=\{f\in\mathcal{E}(U,\CC): \|f\|_{M,U,h}<+\infty\}$ and $\mathcal{E}_{\omega,l}(U,\CC):=\{f\in\mathcal{E}(U,\CC): \|f\|_{\omega,U,l}<+\infty\}$ (Banach spaces). Instead of $\mathcal{E}_{\{M\}}(U,\CC)$ respectively $\mathcal{E}_{\{\omega\}}(U,\CC)$ we will write now $\mathcal{E}_{\{M\}}(U)$ respectively $\mathcal{E}_{\{\omega\}}(U)$ and we have $\mathcal{A}_{\{M\}}(S_{\gamma})=\mathcal{H}(S_{\gamma})\cap\mathcal{E}_{\{M\}}(S_{\gamma})$ and $\mathcal{A}_{\{\omega\}}(S_{\gamma})=\mathcal{H}(S_{\gamma})\cap\mathcal{E}_{\{\omega\}}(S_{\gamma})$, see \cite[Section 2.2]{Thilliezdivision}.\vspace{6pt}

If $M\in\hyperlink{LCset}{\mathcal{LC}}$ and $N$ arbitrary, then $M\hyperlink{precsim}{\precsim}N$ if and only if $\mathcal{E}_{\{M\}}\subseteq\mathcal{E}_{\{N\}}$. For $\sigma,\tau\in\hyperlink{omset1}{\mathcal{W}}$ we get $\sigma\hyperlink{ompreceq}{\preceq}\tau\Leftrightarrow\mathcal{E}_{\{\sigma\}}\subseteq\mathcal{E}_{\{\tau\}}$, see \cite[Corollary 5.17]{compositionpaper}.\vspace{6pt}

In the literature frequently also local classes are considered, here in \eqref{defofm} and \eqref{defofomega} one requires that for each compact set $K\subseteq U$ there exists $h>0$ respectively $l>0$ such that $\|f\|_{M,K,h}<+\infty$ respectively $\|f\|_{\omega,K,l}<+\infty$. Write $\mathcal{E}^{\text{loc}}_{\{M\}}(U)$ respectively $\mathcal{E}^{\text{loc}}_{\{\omega\}}(U)$ for such classes and \eqref{repr1} turns into
\begin{equation*}\label{repr2}
\mathcal{E}^{\text{loc}}_{\{M\}}(U,\CC):=\underset{K\subseteq U}{\varprojlim}\;\underset{h>0}{\varinjlim}\;\mathcal{E}_{M,h}(K,\CC)\hspace{25pt}\mathcal{E}^{\text{loc}}_{\{\omega\}}(U,\CC):=\underset{K\subseteq U}{\varprojlim}\;\underset{l>0}{\varinjlim}\;\mathcal{E}_{\omega,l}(K,\CC).
\end{equation*}

Similarly the Beurling type classes $\mathcal{E}^{\text{loc}}_{(M)}$ and $\mathcal{E}_{(M)}$ respectively $\mathcal{E}^{\text{loc}}_{(\omega)}$ and $\mathcal{E}_{(\omega)}$ can be introduced by replacing $\exists\;h>0$ respectively $\exists\;l>0$ by $\forall\;h>0$ respectively $\forall\;l>0$ and replacing $\underset{h>0}{\varinjlim}$ respectively $\underset{l>0}{\varinjlim}$ by $\underset{h>0}{\varprojlim}$ respectively $\underset{l>0}{\varprojlim}$.

\subsection{Preliminaries}\label{ultraholomorphicextensionprelim}
For the proof of the main theorem we will need to deal with precise estimations of nonzero functions in $\mathcal{E}_{\{\omega\}}$ which are flat at the origin $0\in\RR^d$, i.e. $\mathcal{B}(f)=(f^{(j)}(0))_{j\in\NN}=0$ (and $\omega$ is of course required to be nonquasianalytic, i.e. \hyperlink{omnq}{$(\omega_{\text{nq}})$} is satisfied). More precisely we have to transfer \cite[$(16)$]{Thilliezdivision} to the weight function setting and are proving estimations of $\sup_{j\in\NN^d, s\in\RR^d}|f^{(j)}(s)|$ in terms of the sequences $W^x$ and its associated functions $h_{w^x}$.

This is the aim of this section but in fact we are dealing with a more general situation. On the one hand we consider the local case $\mathcal{E}^{\text{loc}}_{\{M\}}$ respectively $\mathcal{E}^{\text{loc}}_{\{\omega\}}$ and for weight sequences not necessarily having \hyperlink{mg}{$(\text{mg})$} (mixed setting). On the other hand we replace $\{0\}\in\RR^d$ by a general closed subset $X$. More precisely we are generalizing \cite[Lemma 2.9]{thilliez-lojaideals} first to the (mixed) weight sequence case and then to the weight function setting.\vspace{6pt}

To do so we have to introduce some notation. Let $U\subseteq\RR^d$ be nonempty open and $X\subseteq U$ be closed in $U$. Put $\dist(t,X):=\inf\{|t-s|: s\in X\}$ and let $M\in\hyperlink{LCset}{\mathcal{LC}}$ be nonquasianalytic, i.e. \hyperlink{mnq}{$(\on{nq})$} holds true. $I_{X}^{\infty}$ shall denote the ideal of functions in $\mathcal{E}(U)$ which are flat on $X$. Finally we set $I_{X,M}^{\on{loc},\infty}:=I_{X}^{\infty}\cap\mathcal{E}^{\text{loc}}_{\{M\}}(U)$ respectively $I_{X,M}^{\infty}:=I_{X}^{\infty}\cap\mathcal{E}_{\{M\}}(U)$.

Similarly we put $I_{X,\omega}^{\on{loc},\infty}:=I_{X}^{\infty}\cap\mathcal{E}^{\text{loc}}_{\{\omega\}}(U)$ respectively $I_{X,\omega}^{\infty}:=I_{X}^{\infty}\cap\mathcal{E}_{\{\omega\}}(U)$ for any weight $\omega\in\hyperlink{omset}{\mathcal{W}}$ satisfying \hyperlink{omnq}{$(\omega_{\on{nq}})$}.

Note that the Denjoy-Carleman theorem which characterizes nonquasianalyticity for classes of ultradifferentiable functions defined by weight sequences, functions and even matrices (e.g. see \cite[Section 4]{testfunctioncharacterization} and the references therein) holds true for global and local classes simultaneously. We prove the following generalization of \cite[Lemma 2.9]{thilliez-lojaideals} where $M\in\hyperlink{SRset}{\mathcal{SR}}$ has been assumed.

\begin{lemma}\label{Lemma29}
Let $M\in\hyperlink{LCset}{\mathcal{LC}}$ be given and satisfying \hyperlink{mnq}{$(\on{nq})$} (which implies $\lim_{p\rightarrow+\infty}(m_p)^{1/p}=+\infty$). Let $U\subseteq\RR^d$ be nonempty open and $f\in I_{X,M}^{\on{loc},\infty}$. Finally let $N\in\RR_{>0}^{\NN}$ such that $N\ge M$ and \hyperlink{genmg}{$(M,N)_{(\on{mg})}$}. Then for all compact $K\subseteq U$ we get:
\begin{equation}\label{Lemma29equ}
\exists\;H_1,H_2\ge 1\;\forall\;j\in\NN^d\;\forall\;s\in K:\;\;\;|f^{(j)}(s)|\le H_1 H_2^{|j|}N_{|j|}h_n(H_2\dist(s,X)).
\end{equation}
If $f\in I_{X,M}^{\infty}$, then \eqref{Lemma29equ} holds true for any $s\in U$ with uniform constants $H_1$, $H_2$.
\end{lemma}

Note that by the assumption on $N$ the function $h_{n}$ is well-defined and nondecreasing, so we can take w.l.o.g. $H_1=H_2$.

\demo{Proof}
Let $r>0$ be given and put $K_r:=\{w\in U: \dist(w,K)\le r\}$. If $r>0$ is small enough the compact set $K_r$ is contained in $U$. Thus there exist $H,B>0$ such that for all $w\in K_r$ and $j,k\in\NN^d$ we get
\begin{equation}\label{Lemma29equ1}
|f^{(j+k)}(w)|\le H B^{|j+k|}M_{|j+k|}\le H(A B)^{|j|+|k|}N_{|j|}N_{|k|},
\end{equation}
where for the second estimate we have used \hyperlink{genmg}{$(M,N)_{(\on{mg})}$}. Let $s\in K$ and $t\in X$ with $|s-t|=\dist(s,X)=\inf\{|s-w|: w\in X\}$ and distinguish two cases.

If $|s-t|=\dist(s,X)\le r$ then by definition the line segment $[s,t]$ is contained in $K_r$. $f^{(j)}$ is flat at $t$ and so Taylor's formula implies $|f^{(j)}(s)|\le d^{|k|}\sup_{k\in\NN^d, w\in K_r}\frac{|f^{(j+k)}(w)| |s-t|^{|k|}}{|k|!}$. Thus for all $s\in K$ and $j\in\NN^d$ we get
\begin{align*}
|f^{(j)}(s)|&\le d^{|k|}\sup_{k\in\NN^d, w\in K_r}\frac{|f^{(j+k)}(w)| |s-t|^{|k|}}{|k|!}\underbrace{\le}_{\eqref{Lemma29equ1}}d^{|k|}H(A B)^{|j|+|k|}N_{|j|}N_{|k|}\frac{|s-t|^{|k|}}{|k|!}
\\&
\le H(A B)^{|j|}N_{|j|} n_{|k|}(d A B\dist(s,X))^{|k|}.
\end{align*}
We take now the infimum with respect to all $k\in\NN^d$ and get, by definition \eqref{functionhequ1},
$$|f^{(j)}(s)|\le H(A B)^{|j|}N_{|j|} h_{n}(d A B\dist(s,X)),$$
so \eqref{Lemma29equ} follows with $H_1:=H$ and $H_2:=dA B$.
If $|s-t|=\dist(s,X)>r$, then $h_m(B\dist(s,X))\ge h_m(Br)$ holds for any sequence $m$ and $B>0$ and $f\in\mathcal{E}^{\text{loc}}_{\{M\}}(U)$ implies that there exist $H_0,B_0>0$ such that for all $s\in K$ and $j\in\NN^d$ we get
$$|f^{(j)}(s)|\le H_0 B_0^{|j|}M_{|j|}=\frac{H_0}{h_m(B_0r)}B_0^{|j|}M_{|j|}h_m(B_0r)\le\frac{H_0}{h_m(B_0r)}B_0^{|j|}M_{|j|}h_m(B_0\dist(s,X)),$$
i.e. \eqref{Lemma29equ} with $H_1:=\frac{H_0}{h_m(B_0r)}$ and $H_2:=B_0$. Recall that $n\ge m$ implies $h_{n}\ge h_m$, see Lemma \ref{functionhproperties}.

In the global version, if $f\in I_{X,M}^{\infty}$, then the proof is much simpler and reduces to the first case (since \eqref{Lemma29equ1} holds true for all $w\in U$ and uniform constants $H$ and $B$).
\qed\enddemo

\begin{corollary}\label{Lemma29corollary}
Let $\omega\in\hyperlink{omset}{\mathcal{W}}$ be given and satisfying \hyperlink{omnq}{$(\omega_{\on{nq}})$}, let $\Omega=\{W^l: l>0\}$ be the weight matrix associated with $\omega$. Let $U\subseteq\RR^d$ be nonempty open and $f\in I_{X,\omega}^{\on{loc},\infty}$. Then for all compact $K\subseteq U$ we get:
\begin{equation}\label{Lemma29corollaryequ}
\exists\;l>0\;\;\exists\;\widetilde{H}\ge 1\;\forall\;j\in\NN^d\;\forall\;s\in K:\;\;\;|f^{(j)}(s)|\le\widetilde{H}W^{2l}_{|j|}h_{w^{2l}}(\widetilde{H}\dist(s,X)).
\end{equation}

If $f\in I_{X,\omega}^{\infty}$, then the index $l>0$ is valid for the whole set $U$.
\end{corollary}

\demo{Proof}
By \cite[Corollary 4.8]{testfunctioncharacterization} each $W^l$ has \hyperlink{mnq}{$(\on{nq})$}, hence we can apply Lemma \ref{Lemma29} to $M:=W^l$. \eqref{Lemma29equ1} turns then into
\begin{equation*}\label{Lemma29corollaryequ1}
|f^{(j+k)}(w)|\le H W^l_{|j+k|}\le HW^{2l}_{|j|}W^{2l}_{|k|},
\end{equation*}
for some $H\ge 1$, $l>0$ depending on given compact set $K\subseteq U$, all $w\in K_r$ with $r>0$ chosen small enough and all $j,k\in\NN^d$ (taking into account \eqref{newmoderategrowth}). Following then the proof of Lemma \ref{Lemma29} with $A=B=1$ we have shown in the first case \eqref{Lemma29equ} with $N=W^{2l}$, $H_1:=H$ and $H_2:=d$ in the argument of the function $h_{w^{2l}}$ but without a factor $H_2^{|j|}$. In the second case we can put $B_0=1$ and so \eqref{Lemma29equ} follows with $H_1:=\frac{H_0}{h_{w^l}(r)}$ and $H_2:=1$ (again without exponential growth) and since $W^l\le W^{2l}$.

In the global case the index $l>0$ is not depending on the compact set and holds uniformly on whole $U$, see \eqref{repr1}.
\qed\enddemo

{\itshape Remark:} We point out that Lemma \ref{Lemma29} and Corollary \ref{Lemma29corollary} hold true for the corresponding Beurling type classes as well and which have not been considered in \cite[Lemma 2.9]{thilliez-lojaideals}. More precisely, by \cite[Section 4]{testfunctioncharacterization}, these classes are nonquasianalytic, too. We can choose $B$ and $B_0$ in Lemma \ref{Lemma29} as small as desired, i.e. in \eqref{Lemma29equ} we can replace
$\exists\;H_1,H_2\ge 1$ there by for all $H_2>0$ (small) there exists $H_1>0$ (large). Similarly in \eqref{Lemma29corollaryequ} replace $\exists\;l>0\;\;\exists\;\widetilde{H}\ge 1$ by for all $l>0$ (small) there exists $\widetilde{H}\ge 1$ (large).\vspace{6pt}

Since for the proof of our main theorem we are following the (single) weight sequence approach from V. Thilliez in \cite{Thilliezdivision} we will have to apply the (ultradifferentiable) $\mathcal{E}_{\{\omega\}}$-version of the Whitney extension theorem (analogous to \cite[Proposition 1.2.3]{Thilliezdivision} in the weight sequence case). Fortunately the weight function case has been completely characterized in \cite{BonetBraunMeiseTaylorWhitneyextension} for arbitrary compact sets $K$.

In \cite{BonetMeiseTaylorSurjectivity} for $K=\{0\}$ and in \cite{Langenbruch94} for arbitrary compact convex sets an ultradifferentiable extension theorem even in a mixed setting has been characterized. Very recently in a joint work the third author and A. Rainer have succeeded to prove a mixed setting for arbitrary compact sets with a slight restriction on the weights, see \cite{whitneyextensionmixedweightfunction}.

\begin{theorem}\label{Theorem123}
Let $\omega\in\hyperlink{omset1}{\mathcal{W}}$ be a strong weight, i.e. \hyperlink{omsnq}{$(\omega_{\on{snq}})$} holds true. Then we get:
\begin{itemize}
\item[$(i)$] There exists some $d\ge 1$ depending only on the weight $\omega$ and the dimension $n$ of $\RR^n$ such that for any $l>0$ there exists a continuous linear extension operator $E_l$ such that
    $$E_l:\Lambda^n_{\omega,l}\longrightarrow\mathcal{E}_{\omega,dl}(\RR^n)$$
    and so $(\mathcal{B}\circ E_l)(a)=a$ for any $\lambda\in\Lambda^n_{\omega,l}$ with $\mathcal{B}$ denoting here the map $f\mapsto(f^{(j)}(0))_{j\in\NN}$. The extension $E_l(a)$ can be assumed to have compact support which is contained in a prescribed neighborhood of $0$.
\item[$(ii)$] For any bounded open $U\subseteq\RR^n$ with Lipschitz boundary we can find some $d'\ge 1$ depending on the weight $\omega$ and $U$ such that for all $l>0$ there exists a continuous linear extension operator
    $$\widetilde{E}_l:\mathcal{E}_{\omega,l}(U)\longrightarrow\mathcal{E}_{\omega,d'l}(\RR^n)$$
    and such that $\widetilde{E}_l(f)|_U=f$ for any element $f\in\mathcal{E}_{\omega,l}(U)$.
\end{itemize}
\end{theorem}

\demo{Proof}
Use the $\mathcal{E}_{\{\omega\}}$-version of the Whitney extension theorem, for $(i)$ take $K=\{0\}$ and for $(ii)$ take $K=\overline{U}$. Note that the proof of \cite[Theorem 3.9]{BonetBraunMeiseTaylorWhitneyextension} gives the existence of a continuous linear extension operator $E_l$ (by the formula of the extension defining $\widetilde{f}$ on p. $173$) but which has not been mentioned explicitly there.
\qed\enddemo

\subsection{A construction of a convenient ramified weight matrix}\label{ramificationpreparation}
The goal of this Section is to prove a generalization of \cite[Section 5]{sectorialextensions}, in particular of \cite[Theorem 5.3]{sectorialextensions}. More precisely we are working here with an additional ramification parameter $q>0$ which will be needed in case 2 in the proof of our main Theorem \ref{Theorem321} below (and $q=1$ yields the results from \cite[Section 5]{sectorialextensions}).\vspace{6pt}

First we recall \cite[Corollary 5.6]{sectorialextensions} which implies an important fact for strong weight functions.

\begin{lemma}\label{comparisonmaintheorem}
Let $\omega\in\hyperlink{omset0}{\mathcal{W}_0}$ with $\gamma(\omega)>1$ be given and let $\Omega=\{W^x: x>0\}$ be the associated weight matrix. Then $\Omega$ is equivalent with respect to \hyperlink{Mroumapprox}{$\{\approx\}$} to a weight matrix consisting only of strongly log-convex sequences, more precisely to the matrix $\{(p!T^{x,y}_p)_{p\in\NN}: x>0\}$ with $T^{x,y}_p:=\exp(\frac{1}{x}\varphi^{*}_{\omega_{w^y}}(xp))$, $x,y>0$ and $p\in\NN$.

Moreover, in this situation, each $\omega_{w^y}$ belongs to the class \hyperlink{omset}{$\mathcal{W}$}.
\end{lemma}

Let now $\omega\in\hyperlink{omset0}{\mathcal{W}_0}$ with $\gamma(\omega)>1$ be given and $\Omega=\{W^x: x>0\}$, $W^x_p:=\exp(\frac{1}{x}\varphi^{*}_{\omega}(xp))$, shall denote the associated weight matrix. By Lemma \ref{comparisonmaintheorem} we have $\Omega\hyperlink{Mroumapprox}{\{\approx\}}\widehat{\mathcal{S}}$ with $\widehat{\mathcal{S}}:=\{(\widehat{S}_p^x)_{p\in\NN}: x>0\}$, $\widehat{S}^x_p:=p!S^x_p$ and $S^{x}_p:=\exp(\frac{1}{x}\varphi^{*}_{\omega_{w^1}}(xp))$ (i.e. the matrix $\mathcal{S}:=\{S^x: x>0\}$ is associated with the weight $\omega_{w^1}$). We are interested in studying the transformation
\begin{equation}\label{necessarytransformation}
\widehat{S}^x_p\mapsto S^x_p\mapsto(S^x_p)^q\mapsto p!(S^x_p)^q=:\widehat{S}_p^{x,q},\hspace{15pt}q>0.
\end{equation}
For convenience we put in the following $S^{x,q}:=(S^x)^q$, $x,q>0$, and write $\mathcal{S}^q:=\{(S^{x,q}_p)_{p\in\NN}: x>0\}$ respectively $\widehat{\mathcal{S}}^q:=\{(\widehat{S}^{x,q}_p)_{p\in\NN}: x>0\}$. Since $S^x\in\hyperlink{LCset}{\mathcal{LC}}$ for all $x>0$, also $S^{x,q}\in\hyperlink{LCset}{\mathcal{LC}}$ for each $x,q>0$ and moreover $\mathcal{S}$ satisfies both \hyperlink{R-mg}{$(\mathcal{M}_{\{\text{mg}\}})$} and \hyperlink{R-L}{$(\mathcal{M}_{\{\text{L}\}})$}, see Section \ref{weightmatrixfromfunction}. Note that \hyperlink{R-mg}{$(\mathcal{M}_{\{\text{mg}\}})$} and \hyperlink{R-L}{$(\mathcal{M}_{\{\text{L}\}})$} are clearly stable with respect to the mappings from \eqref{necessarytransformation} hence each matrix $\mathcal{S}^q$ and $\widehat{\mathcal{S}}^q$ enjoys both properties too for any $q>0$ (and the same holds true for the Beurling type conditions \hyperlink{B-mg}{$(\mathcal{M}_{(\text{mg})})$} and \hyperlink{R-L}{$(\mathcal{M}_{(\text{L})})$} since $\mathcal{S}$ is associated with $\omega_{w^1}\in\hyperlink{omset}{\mathcal{W}}$, see Section \ref{weightmatrixfromfunction}).\vspace{6pt}

The next result generalizes \cite[Lemma 5.1]{sectorialextensions}, for the convenience of the reader we give the full proof.

\begin{lemma}\label{necessarytransformlemma1}
Let $q>0$ be given, arbitrary but from now on fixed. Then for all $x,y>0$ we get $(\omega^{\iota}_{S^{x,q}})_{\star}\hyperlink{sim}{\sim}(\omega^{\iota}_{S^{y,q}})_{\star}$.
\end{lemma}

\demo{Proof}
First, $(vi)$ in Lemma \ref{assofuncproper} implies that for all $x,y>0$ there exists some $C\ge 1$ such that $-C+C^{-1}\omega_{S^y}(t)\le\omega_{S^x}(t)\le C\omega_{S^y}(t)+C$ for all $t\ge 0$ (since $\omega_{S^x}\hyperlink{sim}{\sim}\omega_{w^1}\hyperlink{sim}{\sim}\omega_{S^y}$ for all $x,y>0$). Then recall \eqref{omegaMspower} which yields $\omega_{S^{x,q}}(t)=\omega_{(S^x)^q}(t)=q\omega_{S^x}(t^{1/q})$ for any $q,x>0$ and $t\ge 0$ and so
\begin{equation*}\label{necessary0}
\forall\;x,y>0\;\exists\;C\ge 1\;\forall\;q>0\;\forall\;t\ge 0:\;\;\;-Cq+C^{-1}\omega_{S^{y,q}}(t)\le\omega_{S^{x,q}}(t)\le C\omega_{S^{y,q}}(t)+Cq,
\end{equation*}
thus by considering $t\mapsto t^{-1}$:
\begin{equation}\label{necessary1}
\forall\;x,y>0\;\exists\;C\ge 1\;\forall\;q>0\;\forall\;t\ge 0:\;\;\;-Cq+C^{-1}\omega_{S^{y,q}}(1/t)\le\omega_{S^{x,q}}(1/t)\le C\omega_{S^{y,q}}(1/t)+Cq.
\end{equation}
Let $x,y>0$ be now arbitrary but fixed and $C\ge 1$ coming from \eqref{necessary1}. By using this inequality we obtain for any $t\ge 0$:
\begin{align*}
(\omega^{\iota}_{S^{x,q}})_{\star}(t)&=\inf_{u>0}\left\{\omega_{S^{x,q}}\left(\frac{1}{u}\right)+ut\right\}\le\inf_{u>0}\left\{C\omega_{S^{y,q}}\left(\frac{1}{u}\right)+ut\right\}+Cq
\\&
=C\inf_{u>0}\left\{\omega_{S^{y,q}}\left(\frac{1}{u}\right)+\frac{ut}{C}\right\}+Cq=C(\omega^{\iota}_{S^{y,q}})_{\star}(t/C)+Cq.
\end{align*}
Proceed similarly for the lower estimate and so, taking into account that each $(\omega^{\iota}_{S^{x,q}})_{\star}$ has \hyperlink{om1}{$(\omega_1)$} and which holds by concavity, we have shown $(\omega^{\iota}_{S^{x,q}})_{\star}\hyperlink{sim}{\sim}(\omega^{\iota}_{S^{y,q}})_{\star}$ for all $x,y,q>0$.
\qed\enddemo

Combining Lemma \ref{necessarytransformlemma1} and \cite[Corollary 3.5]{sectorialextensions} we can show the generalization of \cite[Corollary 5.2]{sectorialextensions}.

\begin{corollary}\label{necessarytransformlemma1corollary}
Let $q>0$ be given, then $\omega_{\widehat{S}^{x,q}}\hyperlink{sim}{\sim}(\omega_{S^{y,q}}^{\iota})_{\star}$ holds for all $x,y>0$ and moreover $\omega_{\widehat{S}^{x,q}}\hyperlink{sim}{\sim}\omega_{\widehat{S}^{y,q}}$ for all $x,y>0$.
\end{corollary}

\demo{Proof}
For any $x>0$ we apply \cite[Corollary 3.5]{sectorialextensions} to the case $m\equiv S^{x,q}\equiv m^{\text{lc}}$, which holds since $S^{x,q}\in\hyperlink{LCset}{\mathcal{LC}}$. Hence we obtain $\omega_{\widehat{S}^{x,q}}\hyperlink{sim}{\sim}(\omega_{S^{x,q}}^{\iota})_{\star}$ and the rest follows from Lemma \ref{necessarytransformlemma1}.
\qed\enddemo

Summarizing everything we are able to prove the following generalization of \cite[Theorem 5.3]{sectorialextensions}.

\begin{theorem}\label{necessarytheorem}
Let $\omega\in\hyperlink{omset0}{\mathcal{W}_0}$ with $\gamma(\omega)>1$ and $q>0$ be given, let $\Omega=\{W^x: x>0\}$ be the associated weight matrix and $\mathcal{S}^q$ and $\widehat{\mathcal{S}}^q$ be the weight matrices as defined above.

Then we have as locally convex vector spaces
\begin{equation}\label{necessarytheoremequ}
\forall\;x>0:\;\;\;\mathcal{E}_{\{\widehat{\mathcal{S}}^q\}}=\mathcal{E}_{\{\omega_{\widehat{S}^{x,q}}\}}=\mathcal{E}_{\{\mathcal{V}^{x,q}\}},
\end{equation}
where $\mathcal{V}^{x,q}$ shall denote the matrix associated with $\omega_{\widehat{S}^{x,q}}$ and the symbol/functor $\mathcal{E}$ can be replaced by $\mathcal{E}^{\on{loc}}$, $\mathcal{A}$ or $\Lambda^n$, $n\in\NN$. In addition $\omega_{\widehat{S}^{x,q}}$ is a strong weight, more precisely
\begin{equation}\label{necessarytheoremequ1}
\forall\;x>0:\;\;\;\gamma(\omega_{\widehat{S}^{x,q}})=q\gamma(\omega)-q+1.
\end{equation}
Note that $\widehat{\mathcal{S}}^q\hyperlink{Mroumapprox}{\{\approx\}}\mathcal{V}^{x,q}$ for all $x>0$ and \eqref{necessarytheoremequ} holds true for the Beurling type spaces as well.
\end{theorem}

\demo{Proof}
Since $\widehat{\mathcal{S}}^q$ satisfies both \hyperlink{R-mg}{$(\mathcal{M}_{\{\text{mg}\}})$} and \hyperlink{R-L}{$(\mathcal{M}_{\{\text{L}\}})$} and since by Corollary \ref{necessarytransformlemma1corollary} all associated functions $\omega_{\widehat{S}^{x,q}}$ (with respect to the parameter $x>0$) are equivalent, we can apply \cite[Corollary 3.17]{testfunctioncharacterization} to obtain the first equality in \eqref{necessarytheoremequ}.

Note that the proof of \cite[Corollary 3.17]{testfunctioncharacterization} shows that $\omega_{\widehat{S}^{x,q}}$ has \hyperlink{om1}{$(\omega_1)$}, i.e. $\omega_{\widehat{S}^{x,q}}\in\hyperlink{omset}{\mathcal{W}}$ and $\gamma(\omega_{\widehat{S}^{x,q}})>0$ by Lemma \ref{omega1weightspace}. So \cite[Theorem 5.14]{compositionpaper} implies the second equality.

In \cite{testfunctioncharacterization} and \cite{compositionpaper} we have considered local classes $\mathcal{E}^{\on{loc}}$, but the proofs are also valid for the symbols $\mathcal{E}$, $\mathcal{A}$ and $\Lambda^n$, $n\in\NN$, by the definition of these classes since we have shown in the proofs the equivalence of the underlying weight structures (see also the explanations in the proof of \cite[Theorem 5.3]{sectorialextensions}). This explains why \eqref{necessarytheoremequ} holds true for the Beurling case as well by taking into account that $\widehat{\mathcal{S}}^q$ satisfies both \hyperlink{B-mg}{$(\mathcal{M}_{(\text{mg})})$} and \hyperlink{R-L}{$(\mathcal{M}_{(\text{L})})$} too since $\mathcal{S}$ does so, see Section \ref{weightmatrixfromfunction}.\vspace{6pt}

To conclude that it is a strong weight we combine \eqref{newindex3}, \eqref{growthindexgammacorollaryequ2} and \eqref{lowerenvelopequ2} to obtain for all $x>0$:
$$\gamma(\omega_{\widehat{S}^{x,q}})=\gamma(\omega_{S^{x,q}})+1=\gamma(\omega_{(S^x)^q})+1=q\gamma(\omega_{S^x})+1=q\gamma(\omega_{w^1})+1=q(\gamma(\omega)-1)+1=q\gamma(\omega)-q+1,$$
and moreover $q\gamma(\omega)-q+1>1\Leftrightarrow\gamma(\omega)>1$.

Finally, \eqref{necessarytheoremequ} applied to $\mathcal{E}^{\on{loc}}$ (or $\mathcal{E}$) implies $\widehat{\mathcal{S}}^q\hyperlink{Mroumapprox}{\{\approx\}}\mathcal{V}^{x,q}$ by using \cite[Proposition 4.6 $(1)$]{compositionpaper} (both matrices are \hyperlink{Msc}{$(\mathcal{M}_{\on{sc}})$} with index set $\mathcal{I}=\RR_{>0}$).
\qed\enddemo

\subsection{The ultraholomorphic extension theorem following V. Thilliez}\label{ultraholomorphicextensioncentralthm}
In this section we state and prove the central theorem of this paper, we transfer \cite[Theorem 3.2.1]{Thilliezdivision} to the weight function setting and are reproving with completely different methods \cite[Theorem 7.7]{sectorialextensions}. For any given $a>0$ we denote by $\lfloor a\rfloor$ its (lower) integer part.

\begin{theorem}\label{Theorem321}
Let $\omega\in\hyperlink{omset0}{\mathcal{W}_0}$ be given with $\gamma(\omega)>1$, i.e. $\omega$ is a strong weight function, and $\Omega=\{W^x: x>0\}$ shall denote the matrix associated with $\omega$. Then for any $\gamma>0$ satisfying $0<\gamma<\gamma(\omega)-1(=\gamma((\omega^{\star})^{\iota})=\gamma(\omega_{W^x})-1=\gamma(\omega_{w^x})$ for any $x>0$, see Corollary \ref{growthindexgammacorollary}) and for all $l>0$ there exists $d''\in\NN_{>0}$ and a continuous linear extension operator
\begin{equation}\label{Theorem321equ0}
E^{\omega}_{\gamma,l}:\Lambda^1_{\omega,l}\longrightarrow\mathcal{A}_{\omega,d''l}(S_{\gamma}),
\end{equation}
i.e. $(\mathcal{B}E^{\omega}_{\gamma,l})(\lambda)=\lambda$ holds true for any sequence $\lambda\in\Lambda^1_{\omega,l}$. Consequently the Borel map $\mathcal{B}:\mathcal{A}_{\{\omega\}}(S_{\gamma})\longrightarrow\Lambda^1_{\{\omega\}}$ is surjective for any $0<\gamma<\gamma(\omega)-1=\gamma((\omega^{\star})^{\iota})$.
\end{theorem}

{\itshape Note:} The proof will show that the number $d''\in\NN_{>0}$ is depending on given index $l$ which is weaker than in \cite[Theorem 3.2.1]{Thilliezdivision} but which seems to be unavoidable in the nonconstant case, see also \cite[Remarque 12, p. 24]{ChaumatChollet94}. It seems that in general \eqref{newexpabsorb} destroys this behavior but which has be used several times, see also the remarks concerning Corollary \ref{Theorem321cor} below.

It suffices to assume $\omega\in\hyperlink{omset0}{\mathcal{W}_0}$ since $\gamma(\omega)>1$ implies \hyperlink{om1}{$(\omega_1)$}, see Lemma \ref{omega1weightspace}.

\demo{Proof}
We follow the proof of \cite[Theorem 3.2.1]{Thilliezdivision} and distinguish between two cases.\vspace{6pt}

{\itshape Case 1:} $\gamma<2$.\vspace{6pt}

{\itshape Part (i) - Construction of ultradifferentiable extensions in $\CC$ with formal holomorphy at $0$.} We identify $\CC$ and $\RR^2$ in the standard way, $\overline{\partial}$ shall denote the Cauchy-Riemann operator, so $\overline{\partial}=\frac{1}{2}\left(\frac{\partial}{\partial x}+i\frac{\partial}{\partial y}\right)$. Let $D$ and $D'$ be two open discs centered at $0$ and such that $\overline{D}\subseteq D'$. For given $l>0$ let $\chi\in\mathcal{E}_{\omega,l}(\CC)$ such that $\supp(\chi)\subseteq D'$ and $\chi(t)=1$ for all $t\in D$. Such a function exists since $\omega$ satisfies \hyperlink{omnq}{$(\omega_{\text{nq}})$}, e.g. see \cite{BraunMeiseTaylor90} and \cite[Section 4]{testfunctioncharacterization}.\vspace{6pt}

Let $\lambda=(\lambda_j)_j\in\Lambda^1_{\omega,l}$ be given and write $\lambda^{\CC}=(\lambda^{\CC}_{j,k})_{(j,k)\in\NN^2}$ for the natural complexification obtained by
$$\sum_{(j,k)\in\NN^2}\lambda^{\CC}_{j,k}\frac{x^jy^k}{j!k!}=\sum_{n\in\NN}\lambda_n\frac{(x+iy)^n}{n!},$$
i.e. $\lambda^{\CC}_{j,k}=i^k\lambda_{j+k}$. So $c:\Lambda^1_{\omega,l}\longrightarrow\Lambda^2_{\omega,l}$, $c:\lambda\mapsto \lambda^{\CC}$, acts as a continuous linear operator with norm~$1$.

Let $E_l:\Lambda^2_{\omega,l}\longrightarrow\mathcal{E}_{\omega,dl}(\CC)$ be the continuous linear extension operator coming from $(i)$ in Theorem \ref{Theorem123} (with constant $d\ge 1$ occurring there) and chosen in such a way that the extension is supported in $D$. Then we put
\begin{equation}\label{gdefinition}
g_{\lambda}:=E_l(\lambda^{\CC}).
\end{equation}

{\itshape Claim:} $\overline{\partial}g_{\lambda}$ is flat at $0$.

First, since $g_{\lambda}=E_l(\lambda^{\CC})$, we get $g^{(j,k)}_{\lambda}(0)=\lambda_{j,k}^{\CC}$ for all $j,k\in\NN$. Hence $$(\overline{\partial}g_{\lambda})^{(j,k)}(0)=\frac{1}{2}\left(g^{(j+1,k)}_{\lambda}(0)+ig^{(j,k+1)}_{\lambda}(0)\right)=\frac{1}{2}\left(\lambda^{\CC}_{j+1,k}+i\lambda^{\CC}_{j,k+1}\right)=\frac{1}{2}\left(i^k\lambda_{j+k+1}+ii^{k+1}\lambda_{j+k+1}\right)=0$$ for all $j,k\in\NN$.\vspace{6pt}

Thus we are able to apply the {\itshape global version} of Lemma \ref{Lemma29} respectively Corollary \ref{Lemma29corollary} for $f:=\overline{\partial}g_{\lambda}$, $U=\CC$ and $X=\{0\}$, so $\dist(s,X)=|s|$. Note that $g_{\lambda}\in\mathcal{E}_{\omega,dl}(\CC)$ by \eqref{gdefinition} and by \eqref{newmoderategrowth} we have $f\in\mathcal{E}_{\omega,2dl}(\CC)$. Hence \eqref{Lemma29corollaryequ} applied for the index $2dl$ implies
\begin{equation}\label{Theorem321equ}
\exists\;C_1\ge 1\;\forall\;j\in\NN^2\;\forall\;\xi\in\CC:\;\;\;|(\overline{\partial}g_{\lambda})^{(j)}(\xi)|\le C_1\|E_l\||\lambda|_{\omega,l}W^{4dl}_{|j|}h_{w^{4dl}}(C_1|\xi|),
\end{equation}
since $\|g_{\lambda}\|_{\omega,\RR^2,2dl}\le \|E_l\||\lambda|_{\omega,l}$ holds true ($\|E_l\|$ denoting the operator norm of $E_l$).\vspace{6pt}

{\itshape Part (ii) - Division of a flat function by a flat ultraholomorphic function.} Let $\tau>0$ and $\xi\in S_{\gamma}$ be given and put
$$H_{a,\tau}(\xi):=G_a(\tau\xi),$$
where $G_a$ denotes the function which we have constructed and studied in Sections \ref{constructionouterfunction} and \ref{sectoriallyflatfunction}. For given $y:=8dl$ (large) we choose the parameter $a=\frac{1}{8y}=\frac{1}{64dl}$ and $x=4y=32dl$ according to the proof of \eqref{Lemma232equ1} in Lemma \ref{Lemma232}. With these choices \eqref{Lemma232equ1} is precisely
\begin{equation}\label{Theorem321equ1}
\exists\;C_2\ge 1\;\forall\;j\in\NN\;\forall\;\xi\in S_{\gamma}:\;\;\;\left|\left(\frac{1}{H_{a,\tau}}\right)^{(j)}(\xi)\right|\le C_2(C_2\tau)^j W^{32dl}_j\frac{1}{h_{w^{8dl}}(C_2^{-1}\tau|\xi|)}.
\end{equation}

Note that the choices of the indices $4dl$, $y$ and $x$ and the constant $C_2$ are not depending on given $\tau$.

We combine now \eqref{Theorem321equ} and \eqref{Theorem321equ1} and estimate for all $(j,k)\in\NN^2$ and $\xi\in S_{\gamma}$ as follows:
\begin{align*}
&\left|\left(\frac{1}{H_{a,\tau}}\overline{\partial}g_{\lambda}\right)^{(j,k)}(\xi)\right|=\left|\sum_{J,K\in\NN^2: J+K=(j,k)}\frac{(j,k)!}{J!K!}\left(\frac{1}{H_{a,\tau}}\right)^{(J)}(\xi)\left(\overline{\partial}g_{\lambda}\right)^{(K)}(\xi)\right|
\\&
\le\sum_{J,K\in\NN^2: J+K=(j,k)}\frac{(j,k)!}{J!K!}C_2(C_2\tau)^{|J|}W^{32dl}_{|J|}\frac{1}{h_{w^{8dl}}(C_2^{-1}\tau|\xi|)}C_1\|E_l\||\lambda|_{\omega,l}W^{4dl}_{|K|}h_{w^{4dl}}(C_1|\xi|)
\\&
\le C_1C_2\|E_l\||\lambda|_{\omega,l}W^{32dl}_{j+k}\frac{h_{w^{4dl}}(C_1|\xi|)}{h_{w^{8dl}}(C_2^{-1}\tau|\xi|)}\sum_{J,K\in\NN^2: J+K=(j,k)}\frac{(j,k)!}{J!K!}(C_2\tau)^{|J|}
\\&
\le C_3\|E\|_l|\lambda|_{\omega,l}(1+\tau C_2)^{j+k}W^{32dl}_{j+k}\frac{h_{w^{4dl}}(C_1|\xi|)}{h_{w^{8dl}}(C_2^{-1}\tau|\xi|)}.
\end{align*}
Here we have used that $W^{4dl}_{|K|}W^{32dl}_{|J|}\le W^{32dl}_{|K|}W^{32dl}_{|J|}\le W^{32dl}_{|J|+|K|}=W^{32dl}_{j+k}$ holds for all $l>0$ and $J,K,(j,k)\in\NN^2$ with $J+K=(j,k)$ and which follows from the log-convexity of each $W^l$.

Let $C_4\ge 1$ denote the universal constant arising in \eqref{functionh2equ1} in Lemma \ref{functionh2}, not depending on $l$, which we will apply to $w^{4dl}$ and $w^{8dl}$ and which explains the choice of the indices. By setting $\tau:=C_1C_2C_4$, which can be done since none of the constants is depending on $\tau$, we get
$$\frac{h_{w^{4dl}}(C_1|\xi|)}{h_{w^{8dl}}(C_2^{-1}\tau|\xi|)}=\frac{h_{w^{4dl}}(C_1|\xi|)}{h_{w^{8dl}}(C_1C_4|\xi|)}\le h_{w^{8dl}}(C_1C_4|\xi|).$$
Finally we put $H_a:=H_{a,C_1C_2C_4}$ and continue the estimate from above:
\begin{align*}
\left|\left(\frac{1}{H_a}\overline{\partial}g_{\lambda}\right)^{(j,k)}(\xi)\right|&\le C_3\|E_l\||\lambda|_{\omega,l}(1+C_1C_2^2C_4)^{j+k}W^{32dl}_{j+k}h_{w^{8dl}}(C_1C_4|\xi|)
\\&
\le C_5\|E_l\||\lambda|_{\omega,l}W^{32A_1dl}_{j+k}h_{w^{8dl}}(C_1C_4|\xi|)\le C_6|\lambda|_{\omega,l}W^z_{j+k}h_{w^{8dl}}(C_6|\xi|)
\\&
\le C_6|\lambda|_{\omega,l}W^z_{j+k}h_{w^z}(C_6|\xi|),
\end{align*}
where we have put $z:=32A_1dl$. We have absorbed again the exponential growth by applying \eqref{newexpabsorb} and note that $A_1\ge 1$ is depending only on $C_1, C_2$ and $C_4$ (via $\tau$). Finally we have used $(i)$ and $(iii)$ from Lemma \ref{functionhproperties}. Note that via \eqref{Theorem321equ1} and the choice of $y$ the constant $A_1$ is now also depending on given index $l>0$ (by constants $C_1$ and $C_2$).\vspace{6pt}

{\itshape Part (iii) - Solution of $\overline{\partial}$-problem.} By the previous estimate of part $(ii)$ above and $h_{w^z}\le 1$ we see that $\frac{1}{H_a}\overline{\partial}g_{\lambda}\in\mathcal{E}_{\omega,z}(S_{\gamma})$ and $\|\frac{1}{H_a}\overline{\partial}g_{\lambda}\|_{\omega,S_{\gamma},z}\le C_6|\lambda|_{\omega,l}$. In the next step we use $(ii)$ of Theorem \ref{Theorem123} for $\Omega:=S_{\gamma}\cap D'$ and put $\nu_{\lambda}:=\widetilde{E}_z(\frac{1}{H_a}\overline{\partial}g_{\lambda})$, so $\nu_{\lambda}\in\mathcal{E}_{\omega,d'z}(\CC)$ and $\|\nu_{\lambda}\|_{\omega,\CC,d'z}\le C_6\|\widetilde{E}_z\||\lambda|_{\omega,l}=C_7|\lambda|_{\omega,l}$. $\nu_{\lambda}$ coincides with $\frac{1}{H_a}\overline{\partial}g_{\lambda}$ on $\Omega$ and $g_{\lambda}$ is assumed to be supported in $D$, hence $\frac{1}{H_a}\overline{\partial}g_{\lambda}$ vanishes on $S_{\gamma}\backslash D$. So it follows that on the whole sector $S_{\gamma}$ and for all $\chi\in\mathcal{E}_{\omega,d'z}(\CC)$, satisfying $\supp(\chi)\subseteq D'$ and $\chi(t)=1$ for all $t\in D$, we get
\begin{equation}\label{Theorem321equ2}
\chi\nu_{\lambda}=\frac{1}{H_a}\overline{\partial}g_{\lambda}.
\end{equation}
The Leibniz rule implies $|(\chi\nu_{\lambda})^{(i)}(\zeta)|\le C_7C_8|\lambda|_{\omega,l}2^{|i|}W^{d'z}_{|i|}$ for all $i\in\NN^2$ and $\zeta\in\CC$ since $W^y_jW^y_k\le W^y_{j+k}$ for all $j,k\in\NN$ and $y>0$ by the log-convexity. Applying once again \eqref{newexpabsorb} we get $|(\chi\nu_{\lambda})^{(i)}(\zeta)|\le C_9|\lambda|_{\omega,l}W^{z_1}_{|i|}$ for all $i\in\NN^2$ (we can take $z_1:=(L(L+1))^{\log(2)}d'z$, $L\ge 1$ the constant arising in \hyperlink{om1}{$(\omega_1)$}, see \eqref{newexpabsorb}) and so
\begin{equation}\label{Theorem321equ3}
\chi\nu_{\lambda}\in\mathcal{E}_{\omega,z_1}(\CC),\hspace{30pt}\|\chi\nu_{\lambda}\|_{\omega,\CC,z_1}\le C_9|\lambda|_{\omega,l}.
\end{equation}
In the next step put $u_{\lambda}:=\mathcal{K}\ast(\chi\nu_{\lambda})$, so the function from above is convolved with the Cauchy kernel $\mathcal{K}(\zeta):=\frac{1}{\pi\zeta}$. We have $\supp(\chi\nu_{\lambda})\subseteq D$, hence $\overline{\partial}(u_{\lambda})=\overline{\partial}(\mathcal{K})\ast(\chi\nu_{\lambda})=\chi\nu_{\lambda}$ and so
$$\forall\;\xi\in\CC:\;\;\;\overline{\partial}(u_{\lambda})(\xi)=\chi\nu_{\lambda}(\xi).$$
Moreover for any $(j,k)\in\NN^2$ and $\xi\in\CC$ we have by changing to polar coordinates
$$|u_{\lambda}^{(j,k)}(\xi)|=|\mathcal{K}\ast(\chi\nu_{\lambda})^{(j,k)}(\xi)|\le\sup_{\zeta\in D}|(\chi\nu_{\lambda})^{(j,k)}(\zeta)|\left|\int_D\frac{1}{\pi\zeta}d\zeta\right|\le 2r\sup_{\zeta\in D}|(\chi\nu_{\lambda})^{(j,k)}(\zeta)|,$$
where $r>0$ denotes the radius of the circle $D$. So $\sup_{\xi\in\CC}|u_{\lambda}^{(j,k)}(\xi)|\le 2r\sup_{\zeta\in D}|(\chi\nu_{\lambda})^{(j,k)}(\zeta)|$ and combining this with \eqref{Theorem321equ3} gives
\begin{equation}\label{Theorem321equ4}
\sup_{\zeta\in\CC}|u_{\lambda}^{(j,k)}(\zeta)|\le 2rC_9|\lambda|_{\omega,l}W^{z_1}_{j+k}=C_{10}|\lambda|_{\omega,l}W^{z_1}_{j+k}.
\end{equation}

{\itshape Part (iv) - Addition of a flat function.}
By \eqref{Lemma232equ}, the definition of $H_a$ and the choice of $\tau$ we get
$$\exists\;C_{11}\ge 1\;\forall\;\xi\in S_{\gamma}\;\forall\;j\in\NN:\;\;\;|H_a^{(j)}(\xi)|\le C_{11}(C_{11}C_1C_2C_4)^jW^{y_1}_jh_{w^{y_1}}(C_{11}|\xi|),$$
where we take $y_1:=4a^{-1}=256dl$. Note that $C_{11}$ is depending on given $l$ since also $a$ is so (see part $(ii)$ above). Together with \eqref{Theorem321equ4} we get now for all $\xi\in S_{\gamma}$ and $(j,k)\in\NN^2$:
\begin{align*}
|(H_a u_{\lambda})^{(j,k)}(\xi)|&\le\sum_{J,K\in\NN^2: J+K=(j,k)}\frac{(j,k)!}{J!K!}\left|H_a^{(J)}(\xi)\right|\left|u^{(K)}_{\lambda}(\xi)\right|
\\&
\le\sum_{J,K\in\NN^2: J+K=(j,k)}\frac{j!k!}{J!K!}C_{10}|\lambda|_{\omega,l}W^{z_1}_{|K|}C_{11}(C_{11}C_1C_2C_4)^{|J|}W^{y_1}_{|J|}h_{w^{y_1}}(C_{11}|\xi|)
\\&
\le C_{10}C_{11}|\lambda|_{\omega,l}W^{y_2}_{j+k}\sum_{J,K\in\NN^2: J+K=(j,k)}\frac{j!k!}{J!K!}(C_{11}C_1C_2C_4)^{|J|}
\\&
\le C_{10}C_{11}|\lambda|_{\omega,l}W^{y_2}_{j+k}(1+C_{11}C_1C_2C_4)^{j+k}\le C_{12}|\lambda|_{\omega,l}W^{A_2y_2}_{j+k},
\end{align*}
where we have put $y_2:=\max\{z_1,y_1\}$ because then $W^{z_1}_{|K|}W^{y_1}_{|J|}\le W^{y_2}_{|K|}W^{y_2}_{|J|}\le W^{y_2}_{|K|+|J|}=W^{y_2}_{j+k}$ holds true by log-convexity. Thus $y_2=\max\{(L(L+1))^{\log(2)}d'z,4a^{-1}\}=\max\{(L(L+1))^{\log(2)}32A_1d'dl,256dl\}$, i.e. $y_2=c_2l$ for some $c_2>0$ (large). In the last step we have absorbed the exponential growth by applying \eqref{newexpabsorb} (with some $A_2\ge 1$ depending only on $C_1, C_2, C_4$ and $C_{11}$ and so on given $l$ through $C_1$, $C_2$ and $C_{11}$). We set $y_3:=A_2y_2$ and so we have shown
\begin{equation}\label{Theorem321equ5}
H_au_{\lambda}\in\mathcal{E}_{\omega,y_3}(S_{\gamma}),\hspace{30pt}\|H_au_{\lambda}\|_{\omega,S_{\gamma},y_3}\le C_{12}|\lambda|_{\omega,l}.
\end{equation}
Put $f_{\lambda}:=g_{\lambda}-H_au_{\lambda}$ and then $f_{\lambda}$ is well-defined and holomorphic in $S_{\gamma}$: We have $\overline{\partial}(f_{\lambda})=\overline{\partial}(g_{\lambda})-\overline{\partial}(H_a)u_{\lambda}-H_a\overline{\partial}(u_{\lambda})=\overline{\partial}(g_{\lambda})-\overline{\partial}(u_{\lambda})=0$, because by \eqref{Theorem321equ2} it follows that $H_a\overline{\partial}(u_{\lambda})=H_a\chi\nu_{\lambda}=H_a\frac{1}{H_a}\overline{\partial}(g_{\lambda})$.\vspace{6pt}

All considered functions $g_{\lambda}$, $\nu_{\lambda}$ and $u_{\lambda}$ depend linearly on $\lambda$ since they were defined in terms of the linear extension operators $E_l$ and $\widetilde{E}_{z}=\widetilde{E}_{32A_1dl}$. \eqref{Theorem321equ5} and the very definition of $g_{\lambda}$ show that $\lambda\mapsto f_{\lambda}$, $\Lambda^1_{\omega,l}\longrightarrow\mathcal{A}_{\omega,y_3}(S_{\gamma})$ is continuous. Note that $y_3=A_2y_2=A_2c_2l$ and $d'':=\frac{y_3}{l}=A_2c_2$ depends on $\gamma$ of the given sector $S_{\gamma}$, the weight $\omega$ and on chosen index $l>0$ as well.

$H_a$ and $H_au_{\lambda}$ are flat at $0$, hence $f_{\lambda}^{(j)}(0)=\frac{\partial^jf_{\lambda}}{\partial x^j}(0)=\frac{\partial^jg_{\lambda}}{\partial x^j}(0)=\lambda_j$ for all $j\in\NN$ and so \eqref{Theorem321equ0} is shown with $d''=A_2c_2$ (without loss of generality we can take $d''\in\NN_{>0}$) and by putting $E^{\omega}_{\gamma,l}(\lambda):=f_{\lambda}$.\vspace{12pt}

{\itshape Case 2:} $\gamma\ge 2$.\vspace{6pt}

We follow the ramification arguments of \cite[Theorem 3.2.1]{Thilliezdivision} and to do so we use the weight structures $\mathcal{S}^{1/q}$ and $\widehat{\mathcal{S}}^{1/q}$ from Section \ref{ramificationpreparation} and apply Theorem \ref{necessarytheorem}. More precisely we apply it to the choice $q:=\lfloor\frac{\gamma}{2}\rfloor+1\ge 2$.

Since $\mathcal{S}$ is associated with the weight $\omega_{w^1}\in\hyperlink{omset}{\mathcal{W}}$, we obtain for each $l>0$ and $j,k\in\NN$ that $S^l_jS^l_k\le S^l_{j+k}\le S^{2l}_jS^{2l}_k$, hence by iterating $(S^l_j)^q\le S^l_{qj}\le(S^{2^{q-1}l}_j)^q$ and so
\begin{equation}\label{Theorem321equ6}
\forall\;l>0\;\forall\;j\in\NN:\;\;\;S^l_j\le(S^l_{qj})^{1/q}\le S^{2^{q-1}l}_j.
\end{equation}

Let $\lambda\in\Lambda^1_{\omega,l}$ be given for some $l>0$ arbitrary but from now on fixed and consider $\lambda^{\ast}$ defined by
\begin{equation}\label{Theorem321equ7}
\lambda^{\ast}_{qj}:=\lambda_j\frac{(qj)!}{j!},\hspace{30pt}\lambda^{\ast}_{qj+k}=0\;\;\;\forall\;j\in\NN, k=1,\dots,q-1.
\end{equation}
By applying the equivalence $\Omega\hyperlink{Mroumapprox}{\{\approx\}}\widehat{\mathcal{S}}$ we get $|\lambda_j|\le CW^l_j\le C h^j\widehat{S}^{l_1}_j$ for some $C,h,l_1>0$ and all $j\in\NN$. Hence, by the first part of \eqref{Theorem321equ6}, we get
\begin{align*}
|\lambda^{\ast}_{qj}|&=|\lambda_j|\frac{(qj)!}{j!}\le C h^j(qj)!S^{l_1}_j\le Ch^j(qj)!(S^{l_1}_{qj})^{1/q}
\\&
=Ch^j(qj)!S^{l_1,1/q}_{qj}=Ch^j\widehat{S}^{l_1,1/q}_{qj}
\end{align*}
and which proves $\lambda^{\ast}\in\Lambda^1_{\widehat{S}^{l_1,1/q},h^{1/q}}$, i.e. $\lambda^{\ast}\in\Lambda^1_{\{\widehat{\mathcal{S}}^{1/q}\}}$. By \eqref{necessarytheoremequ} we obtain now $\lambda^{\ast}\in\Lambda^1_{\{\omega_{\widehat{S}^{x,1/q}}\}}$ for any $x>0$.

For convenience we put in the following $\tau:=\omega_{\widehat{S}^{1,q}}$, moreover let $\mathcal{V}:=\{V^y: y>0\}$ denote the weight matrix associated with $\tau$, hence $\widehat{\mathcal{S}}^{1/q}\hyperlink{Mroumapprox}{\{\approx\}}\mathcal{V}$ as shown in Theorem \ref{necessarytheorem}. So $\lambda^{\ast}\in\Lambda^1_{\tau,l_2}$ holds true for some $l_2>0$ since the equivalence \eqref{necessarytheoremequ} implies $h^j\widehat{S}^{l_1,1/q}_{qj}\le C_1 V^{l_2}_{qj}$. Gathering everything it shows for all $j\in\NN$
$$\frac{|\lambda^{\ast}_{qj}|}{V^{l_2}_{qj}}\le C_1\frac{|\lambda^{\ast}_{qj}|}{h^j(qj)!S^{l_1,1/q}_{qj}}\le C_1\frac{|\lambda^{\ast}_{qj}|}{h^j(qj)!S^{l_1}_j}=C_1\frac{|\lambda_j|}{h^j\widehat{S}^{l_1}_j}\le C_1\frac{|\lambda_j|}{W^l_j},$$
i.e. $|\lambda^{\ast}|_{\tau,l_2}\le C_1|\lambda|_{\omega,l}$ and which proves the continuity of $\lambda\mapsto\lambda^{\ast}$.\vspace{6pt}

The aim is to apply the first case to the weight $\tau$. By Theorem \ref{necessarytheorem}, $\gamma(\omega)>1$ implies $\gamma(\tau)>1$ and by \eqref{necessarytheoremequ1} we obtain $\frac{\gamma}{q}<\gamma(\tau)-1=\frac{1}{q}\gamma(\omega)-\frac{1}{q}\Leftrightarrow\gamma<\gamma(\omega)-1=\gamma((\omega^{\star})^{\iota})$ and which holds true for all values under consideration $2\le\gamma<\gamma((\omega^{\star})^{\iota})$.

Then replace in {\itshape case 1} above $l$ by $l_2$, $\omega$ by $\tau$ and $\gamma$ by $\gamma':=\frac{\gamma}{q}<2$ (which holds by the choice of $q$ above), i.e. the new sector $S_{\gamma'}$ is a subset of $\CC$.\vspace{6pt}

Hence we obtain an extension $h_{\lambda}\in\mathcal{A}_{\tau,l_3}(S_{\gamma'})$ with $l_3=d''l_2$, $d''\in\NN_{>0}$ also coming from the first part of the proof and depending on $l_2$. Since $\lambda\mapsto\lambda^{\ast}$ is linear and continuous the map $\lambda\mapsto h_{\lambda}$, $\Lambda^1_{\omega,l}\rightarrow\mathcal{A}_{\tau,l_3}(S_{\gamma'})$, depends linearly and continuously on given $\lambda\in\Lambda^1_{\omega,l}$ and satisfies
\begin{equation}\label{Theorem321equ8}
\mathcal{B}(h_{\lambda})=\lambda^{\ast}.
\end{equation}
Let $\xi\in S_{\gamma}$ and put $f_{\lambda}(\xi):=h_{\lambda}(\xi^{1/q})$ which is holomorphic in $S_{\gamma}$ and bounded by the supremum norm of $h_{\lambda}$. Consider the differential operator $Y:=q^{-1}\xi^{1-q}\frac{\partial}{\partial\xi}$ on $\CC\backslash\{0\}$. By the chain rule we get
\begin{equation}\label{Theorem321equ9}
\forall\;\xi\in S_{\gamma'}\;\forall\;j\in\NN_{>0}:\;\;\;f_{\lambda}^{(j)}(\xi^q)=(Y^jh_{\lambda})(\xi).
\end{equation}
As pointed out in \cite[Theorem 3.2.1]{Thilliezdivision} we have estimates for the coefficients of $Y^j=\sum_{k=1}^jY_{j,k}(\xi)\frac{\partial^k}{\partial\xi^k}$ as follows:
\begin{equation}\label{Theorem321equ10}
\forall\;\xi\;\forall\;k,j\in\NN_{>0},\;1\le k\le j:\;\;\;|Y_{j,k}(\xi)|\le(4q^{-1})^j2^{j-k}(j-k)!|\xi|^{k-qj}.
\end{equation}
For the following considerations we define the polynomial $P_{\lambda,j}$ of degree at most $q(j-1)$ by
\begin{equation}\label{Theorem321equ11}
P_{\lambda,j}(\xi):=\sum_{i=0}^{j-1}\lambda_i\frac{\xi^{qi}}{i!}.
\end{equation}
By \eqref{Theorem321equ7} and \eqref{Theorem321equ8} we obtain
\begin{align*}
P_{\lambda,j}(\xi)&=\sum_{i=0}^{j-1}\lambda_i\frac{\xi^{qi}}{i!}=\sum_{i=0}^{j-1}\frac{\lambda^{\ast}_{qi}i!}{(qi)!}\frac{\xi^{qi}}{i!}=\sum_{i=0}^{j-1}\lambda^{\ast}_{qi}\frac{\xi^{qi}}{(qi)!}=\sum_{p=0}^{q(j-1)}\lambda^{\ast}_p\frac{\xi^p}{p!}
\\&
=\sum_{p=0}^{q(j-1)}h_{\lambda}^{(p)}(0)\frac{\xi^p}{p!}=\sum_{p=0}^{qj-1}h_{\lambda}^{(p)}(0)\frac{\xi^p}{p!}.
\end{align*}
The last equality holds by the definition of $\lambda^{\ast}$. In the next step the Taylor formula applied for $h_{\lambda}$ between $0$ and any $\xi\in S_{\gamma'}$ gives, for all $1\le k\le j$,
\begin{align*}
\left|\frac{\partial^k}{\partial \xi^k}(h_{\lambda}-P_{\lambda,j})(\xi)\right|&\le\sup_{\zeta\in(0,\xi)}|h^{(qj)}_{\lambda}(\zeta)|\frac{|\xi|^{qj-k}}{(qj-k)!}\le C_2|\lambda|_{\omega,l}V^{l_3}_{qj}\frac{|\xi|^{qj-k}}{(qj-k)!}\le C_3b^{qj}|\lambda|_{\omega,l}\widehat{S}^{l_4,1/q}_{qj}\frac{|\xi|^{qj-k}}{(qj-k)!}
\\&
=C_3b^{qj}|\lambda|_{\omega,l}(qj)!(S^{l_4}_{qj})^{1/q}\frac{|\xi|^{qj-k}}{(qj-k)!}\le  C_3b^{qj}|\lambda|_{\omega,l}2^{qj}k!(S^{l_4}_{qj})^{1/q}|\xi|^{qj-k}
\\&
\le C_3|\lambda|_{\omega,l}(2b)^{qj}k!S^{2^{q-1}l_4}_j|\xi|^{qj-k},
\end{align*}
where in the third estimate we have used $\widehat{\mathcal{S}}^{1/q}\hyperlink{Mroumapprox}{\{\approx\}}\mathcal{V}$ (see Section \ref{ramificationpreparation}) and for the last one the second part of \eqref{Theorem321equ6}. This and \eqref{Theorem321equ10} imply for any $\xi\in S_{\gamma'}$ and $j\in\NN$:
\begin{align*}
&\left|Y^j(h_{\lambda}-P_{\lambda,j})(\xi)\right|
\\&
\le\sum_{k=1}^j\left|Y_{j,k}(\xi)\right|\left|\frac{\partial^k}{\partial \xi^k}(h_{\lambda}-P_{\lambda,j})(\xi)\right|\le(4q^{-1})^j\sum_{k=1}^j2^{j-k}(j-k)!|\xi|^{k-qj}C_3|\lambda|_{\omega,l}(2b)^{qj}k!S^{2^{q-1}l_4}_j|\xi|^{qj-k}
\\&
\le C_3(8q^{-1})^j(2b)^{qj}|\lambda|_{\omega,l}j!S^{2^{q-1}l_4}_j\sum_{k=1}^j2^{-k}\le C_3(8q^{-1}(2b)^q)^j|\lambda|_{\omega,l}\widehat{S}^{2^{q-1}l_4}_j
\\&
\le C_4|\lambda|_{\omega,l}(8q^{-1}(2bb_1)^q)^j|\lambda|_{\omega,l}W^{l_5}_j\le C_5|\lambda|_{\omega,l}W^{l_6}_j,
\end{align*}
where in the second last estimate we have used $\Omega\hyperlink{Mroumapprox}{\{\approx\}}\widehat{\mathcal{S}}$ and finally \eqref{newexpabsorb}. \eqref{Theorem321equ11} implies $P_{\lambda,j}(\xi)=Q_{\lambda,j}(\xi^q)$ where the polynomial $Q_{\lambda,j}$ has degree at most $j-1$, so
\begin{equation}\label{Theorem321equ12}
(Y^jP_{\lambda,j})(\xi)=Q_{\lambda,j}^{(j)}(\xi^q)=0.
\end{equation}
\eqref{Theorem321equ9}, \eqref{Theorem321equ12} and the previous estimate show for any $\xi\in S_{\gamma'}$ and all $j\in\NN$:
\begin{align*}
|f^{(j)}_{\lambda}(\xi^q)|&=|(Y^jh_{\lambda})(\xi)|\le|Y^j(h_{\lambda}-P_{\lambda,j})(\xi)|+|(Y^jP_{\lambda,j})(\xi)|\le  C_5|\lambda|_{\omega,l}W^{l_6}_j,
\end{align*}
i.e. we have shown $f_{\lambda}\in\mathcal{A}_{\omega,l_6}(S_{\gamma})$ and $\lambda\mapsto f_{\lambda}$ is linear and continuous. Finally \eqref{Theorem321equ8} implies that the Taylor series of $h_{\lambda}$ at $0$ is given by $\sum_{j\in\NN}\lambda_j\frac{\xi^{qj}}{j!}$, hence $f_{\lambda}$ has the expansion $\sum_{j\in\NN}\lambda_j\frac{\xi^j}{j!}$ (by replacing $\xi\mapsto\xi^{1/q}$) and so $\mathcal{B}(f_{\lambda})=\lambda$ as desired.
\qed\enddemo

Using this result we can reprove \cite[Theorem 3.2.1]{Thilliezdivision} independently by using the associated weight function.

\begin{corollary}\label{Theorem321cor}
Let $\widehat{M}=(p!M_p)_{p\in\NN}$ be given such that $\widehat{M}\in\hyperlink{LCset}{\mathcal{LC}}$, $\gamma(\omega_{\widehat{M}})>1$ and satisfying \hyperlink{mg}{$(\on{mg})$}. Then for any $\gamma>0$ satisfying $0<\gamma<\gamma(\omega_{\widehat{M}})-1=\gamma(\omega_M)=\gamma(M)=\gamma(\widehat{M})-1$ and for all $h>0$ there exists $h_1>0$ and a continuous linear extension operator
\begin{equation}\label{Theorem321corequ0}
E^{\widehat{M}}_{\gamma,h}:\Lambda^1_{\widehat{M},h}\longrightarrow\mathcal{A}_{\widehat{M},h_1}(S_{\gamma}),
\end{equation}
i.e. $(\mathcal{B}E^{\widehat{M}}_{\gamma,h})(\lambda)=\lambda$ holds true for any sequence $\lambda\in\Lambda^1_{\widehat{M},h}$. Consequently the Borel map $\mathcal{B}:\mathcal{A}_{\{\widehat{M}\}}(S_{\gamma})\longrightarrow\Lambda^1_{\{\widehat{M}\}}$ is surjective for any $0<\gamma<\gamma(\omega)-1=\gamma(M)=\gamma(\widehat{M})-1$.
\end{corollary}

\demo{Proof}
By assumption we have $\omega_{\widehat{M}}\in\hyperlink{omset1}{\mathcal{W}}$ and by Theorem \ref{Theorem321} for any $\gamma>0$ satisfying $0<\gamma<\gamma((\omega_{\widehat{M}}^{\star})^{\iota})$ the Borel map $\mathcal{B}:\mathcal{A}_{\{\omega_{\widehat{M}}\}}(S_{\gamma})\longrightarrow\Lambda^1_{\{\omega_{\widehat{M}}\}}$ is surjective. Denote by $\Omega=\{W^x: x>0\}$ the weight matrix associated with $\omega_{\widehat{M}}$. By $(iii)$ in Lemma \ref{assofuncproper}, $(b)$ in Remark \ref{importantremark} and \eqref{omegaMassofctmatrix} we see $\widehat{M}\equiv W^1\hyperlink{approx}{\approx}W^x$ for all $x>0$, in particular $\Omega$ is constant. Thus $\mathcal{A}_{\{\omega_{\widehat{M}}\}}(S_{\gamma})=\mathcal{A}_{\{\Omega\}}(S_{\gamma})=\mathcal{A}_{\{\widehat{M}\}}(S_{\gamma})$ for any $\gamma>0$ and analogously for the class $\Lambda^1$.

More precisely, let now $\lambda\in\Lambda^1_{\widehat{M},h}\equiv\Lambda^1_{W^1,h}$ be given for some $h>0$ arbitrary (large) but fixed. First we apply \eqref{newexpabsorb} to get $h^jW^1_j\le DW^x_j$ for some $D,x>0$ and all $j\in\NN$. Then we have the extension
$E^{\omega_{\widehat{M}}}_{\gamma,x}:\Lambda^1_{W^x,1}\equiv\Lambda^1_{\omega_{\widehat{M}},x}\longrightarrow\mathcal{A}_{\omega_{\widehat{M}},d''x}(S_{\gamma})\equiv\mathcal{A}_{W^{d''x},1}(S_{\gamma})$ for some $d''$ depending on $x$ (and so on $h$). Finally $W^{xd''}_j\le D_1h_1^jW^1=h_1^j\widehat{M}_j$ for some $D_1,h_1>0$ and all $j\in\NN$ since $\Omega$ is constant.

Thus we have shown the existence of an extension operator $E^{\widehat{M}}_{\gamma,h}:\Lambda^1_{\widehat{M},h}\longrightarrow\mathcal{A}_{\widehat{M},h_1}(S_{\gamma})$, i.e. \eqref{Theorem321corequ0}.
\qed\enddemo

Let us point out that if $\widehat{M}\in\hyperlink{LCset}{\mathcal{LC}}$ is satisfying \hyperlink{mg}{$(\on{mg})$}, then $\gamma(\omega_{\widehat{M}})>1$ is equivalent to \hyperlink{gamma1}{$(\gamma_1)$} for $\widehat{M}$ as we show in \cite[Section 4]{firstindexpaper}. Moreover in \cite[Corollary 4.8]{firstindexpaper} we prove $\gamma(\widehat{M})=\gamma(\omega_{\widehat{M}})$ and so $\gamma(M)=\gamma(\omega_M)=\gamma(\omega_{\widehat{M}})-1$ (hence the upper bound of the admissible opening of the sectors is the same as in \cite[Theorem 3.2.1]{Thilliezdivision}).

By \cite{petzsche} for any $\widehat{M}$ as considered in Corollary \ref{Theorem321cor} there exists a sequence $N$ such that $N\hyperlink{approx}{\approx}\widehat{M}$ and $N\in\hyperlink{SRset}{\mathcal{SR}}$.\vspace{6pt}

A disadvantage of this result is that one can not conclude $h_1=ch$ for some $c$ not depending on given $h>0$ as in \cite[Theorem 3.2.1]{Thilliezdivision}. For this note that $x=ch$ will not be possible in general for some universal $c$. For given $h$, as mentioned in \eqref{newexpabsorb}, let $a\in\NN$ be chosen minimal to have $\exp(a)\ge h$ and we have to consider the choice $x=L_1^a$, $L_1=L(L+1)\ge L$ and $L\ge 1$ denoting the constant arising in \hyperlink{om1}{$(\omega_1)$}. To guarantee $ch\ge L_1^a$ we would need $\log(c)/\log(h)+1\ge\log(L_1)$ and which leads to $L_1\le\exp(1)$ as $h\rightarrow\infty$ since $c$ should be not depending on $h$. Hence $L\le\sqrt{e}$ would follow and which is a contradiction in general since $\gamma(\omega_{\widehat{M}})>1$ and \eqref{omegasnqiteration} only imply $\omega(2t)\le 2\omega(t)+C$ for all $t\ge 0$ and some $C>0$.\vspace{6pt}

Moreover we do not see that $d''$ is not depending on given $x$. But even if so we would require $h_1\le c_1 h$ for some universal constant $c_1$. Inspecting the proof of \cite[$(5.11)$]{compositionpaper} (respectively \cite[$(3.4.1)$]{dissertation}) yields that we have to iterate the last estimate there $j$-times such that $2^j\le d''x<2^{j+1}$. Then we would require $h_1=H^{j+1}\le c_1h$ for some universal constant $H$ (large), but $\log(c_1)+\log(h)\ge(j+1)\log(H)\ge(j+1)\log(2)\ge\log(x)+\log(d'')=a\log(L)+\log(d'')\ge\log(h)\log(L_1)+\log(d'')$ again leads to a contraction as $h\rightarrow+\infty$.\vspace{6pt}

Corollary \ref{Theorem321cor} should also be compared with \cite[Remark 7.9]{sectorialextensions} where we have been able to omit this problem by using the first main extension result \cite[Theorem 7.4]{sectorialextensions} and applying complex methods for the proof there.

\section{The Beurling case}\label{Beurlingcase}
The aim of this section is to prove Theorem \ref{Theorem321} for the Beurling case as well and the strategy is to reduce this situation to the Roumieu case as it has been indicated in \cite[Section 3.4]{Thilliezdivision} for the single weight sequence case. For this we generalize \cite[Lemma 4.4]{BonetBraunMeiseTaylorWhitneyextension} as follows.

\begin{lemma}\label{lemma44BBMTWhitney}
Let $\omega\in\hyperlink{omset0}{\mathcal{W}_0}$ with $\gamma(\omega)>1$. Let $h:[0,+\infty)\rightarrow[0,+\infty)$ be arbitrary and satisfying $\omega(t)=o(h(t))$ as $t\rightarrow+\infty$. Then for all $0<\gamma<\gamma(\omega)$ there exists a function $\sigma_{\gamma}=\sigma\in\hyperlink{omset0}{\mathcal{W}_0}$ satisfying $\omega(t)=o(\sigma(t))$ and $\sigma(t)=o(h(t))$ as $t\rightarrow+\infty$ and finally such that $\gamma(\sigma)>\gamma$.
\end{lemma}

\demo{Proof}
Let $0<\gamma<\gamma(\omega)$ be given, then by definition \eqref{newindex2} there exists some $K>1$ such that $\limsup_{t\rightarrow+\infty}\frac{\omega(K^{\gamma}t)}{\omega(t)}<K$ holds true.

Since $\omega$ tends to infinity, respectively by \hyperlink{om3}{$(\omega_3)$}, we see that $\lim_{t\rightarrow+\infty}h(t)=+\infty$. We are going to introduce now inductively a (strictly increasing) sequence $(x_n)_{n\ge 1}$ with $x_1=0$, $\omega(x_2)>0$ (so $x_2>1$ by normalization of $\omega$), and increasing fast enough such that the following conditions are satisfied:
\begin{equation}\label{lemma44BBMTWhitneyequ2}
\forall\;n\ge 1:\;\;\;x_{n+1}\ge K^{\gamma}x_n,
\end{equation}
\begin{equation}\label{lemma44BBMTWhitneyequ3}
\omega(x_{n+1})\ge 2^{n+1-i}\omega(x_i),\;\;\;\forall\;1\le i\le n,
\end{equation}
\begin{equation}\label{lemma44BBMTWhitneyequ4}
\forall\;n\ge 1\;\forall\;x\ge x_n:\;\;\;h(x)\ge n^2\omega(x).
\end{equation}
Now define the weight $\sigma$ by
\begin{equation}\label{lemma44BBMTWhitneyequ}
\sigma(x):=n\omega(x)-\sum_{i=1}^n\omega(x_i),\hspace{15pt}\text{for}\;x\in[x_n,x_{n+1}),\;\;n\ge 1.
\end{equation}
Hence $\sigma$ is continuous, nondecreasing, tending to infinity and has \hyperlink{om4}{$(\omega_4)$}. If $x\in[x_1,x_2)=[0,x_2)$, then $\sigma(x)=\omega(x)-\omega(x_1)=\omega(x)$ and so $\sigma$ is normalized, too.

If $x\in[x_n,x_{n+1})$, $n\ge 2$, then by applying \eqref{lemma44BBMTWhitneyequ3} we get $\frac{\omega(x_i)}{\omega(x)}\le\frac{\omega(x_i)}{\omega(x_n)}\le 2^{i-n}$ for all $1\le i\le n-1$. Moreover $\sum_{i=1}^{n-1}2^{i-n}=2^{-n}\sum_{i=1}^{n-1}2^i=2^{-n}(2^n-2)=1-2^{-n+1}\le 1$ and so
\begin{equation}\label{lemma44BBMTWhitneyequ5}
\sigma(x)=\left(n-\sum_{i=1}^n\frac{\omega(x_i)}{\omega(x)}\right)\omega(x)\ge\left(n-\sum_{i=1}^{n-1}2^{i-n}-1\right)\omega(x)\ge(n-2)\omega(x),
\end{equation}
which implies
\begin{equation*}\label{lemma44BBMTWhitneyequ0}
\forall\;n\ge 3\;\forall\;x\in[x_n,x_{n+1}):\;\;\;\omega(x)\le\frac{1}{n-2}\sigma(x).
\end{equation*}
By this we get $\omega(t)=o(\sigma(t))$ as $t\rightarrow+\infty$ and $\sigma$ has \hyperlink{om3}{$(\omega_3)$} too and altogether $\sigma\in\hyperlink{omset0}{\mathcal{W}_0}$. On the other hand, \eqref{lemma44BBMTWhitneyequ} and \eqref{lemma44BBMTWhitneyequ4} imply
$$\forall\;n\ge 1\;\forall\;x\in[x_n,x_{n+1}):\;\;\;\sigma(x)\le n\omega(x)\le\frac{1}{n}h(x),$$
and consequently $\sigma(t)=o(h(t))$ as $t\rightarrow+\infty$.\vspace{6pt}


Let $\varepsilon>0$ be given, then there exists $N=N(\varepsilon)\in\NN$ large ($N\ge 3$) such that
\begin{equation}\label{lemma44BBMTWhitneyequ6}
\forall\;t\ge x_N:\;\;\;\omega(K^{\gamma}t)\le(K-\varepsilon)\omega(t).
\end{equation}

For all $x_{N} \leq x_n\leq t<x_{n+1}$ we see that $K^\gamma t\geq K^\gamma x_n \geq x_n$ and, by \eqref{lemma44BBMTWhitneyequ2}, $K^\gamma t<K^\gamma x_{n+1}\leq x_{n+2}$. This means that $K^\gamma t$ either belongs to $[x_n,x_{n+1})$ or to $[x_{n+1},x_{n+2})$. Then by~\eqref{lemma44BBMTWhitneyequ}, ~\eqref{lemma44BBMTWhitneyequ6} and \eqref{lemma44BBMTWhitneyequ5} we observe that
$$\sigma(K^\gamma t)\leq (n+1) \omega (K^\gamma t)\leq (n+1) (K-\varepsilon)\omega(t) \leq \frac{n+1}{n-2} (K-\varepsilon) \sigma(t).$$
Since $\frac{n+1}{n-2}\rightarrow 1$ as $n\rightarrow+\infty$ we have verified $\limsup_{t\rightarrow+\infty}\frac{\sigma(K^{\gamma}t)}{\sigma(t)}\le K-\varepsilon<K$, i.e. $(P_{\sigma,\gamma})$ holds true with constant $K$ and so $\gamma(\sigma)>\gamma$ (since requirement $(P_{\sigma,\gamma})$ is an open condition as seen in Section \ref{growthindexgamma}).
\qed\enddemo



\begin{theorem}\label{Theorem321Beurling}
Let $\omega\in\hyperlink{omset0}{\mathcal{W}_0}$ be given with $\gamma(\omega)>1$, i.e. $\omega$ is a strong weight function. Then for any $\gamma>0$ satisfying $0<\gamma<\gamma(\omega)-1(=\gamma((\omega^{\star})^{\iota}))$ the Borel map $\mathcal{B}:\mathcal{A}_{(\omega)}(S_{\gamma})\longrightarrow\Lambda^1_{(\omega)}$ is surjective.
\end{theorem}

\demo{Proof}
Let $\lambda=(a_p)_{p\in\NN}\in\Lambda^1_{(\omega)}$ be given. Now we follow the lines of the proof of \cite[Theorem 4.5]{BonetBraunMeiseTaylorWhitneyextension} for $\lambda$ instead of considering a Whitney jet $F$, with $K=\{0\}$ and where we put $g(t):=\log\max\{1,|a_p|\}$ for $p\le t<p+1$, $p\in\NN$. So we construct $h_{\lambda}:[0,+\infty)\rightarrow[0,+\infty)$ such that $\omega(t)=o(h(t))$ as $t\rightarrow+\infty$. We fix $0<\gamma<\gamma(\omega)-1$ and denote by $\sigma$ the weight constructed in Lemma \ref{lemma44BBMTWhitney} for $h_{\lambda}$ and $\gamma+1$. We observe that $\gamma(\sigma)>\gamma+1>1$ then $\gamma(\sigma)>1$ and so the weight $\sigma$ is strong too and $\sigma^{\star}$ is well-defined.

Thus we are able to apply the main result Theorem \ref{Theorem321} to $\sigma$ and $\gamma$, hence the Borel map $\mathcal{B}:\mathcal{A}_{\{\sigma\}}(S_{\gamma})\longrightarrow\Lambda^1_{\{\sigma\}}$ is surjective.

Finally, since $\sigma=o(h_{\lambda}(t))$ as $t\rightarrow+\infty$ we have that $\lambda\in\Lambda^1_{\{\sigma\}}$ and since $\omega(t)=o(\sigma(t))$ as $t\rightarrow+\infty$, by using \cite[Lem. 5.16 (2), Cor. 5.17 (2)]{compositionpaper} we get $\mathcal{A}_{\{\sigma\}}(S_{\gamma})\subseteq\mathcal{A}_{(\omega)}(S_{\gamma})$ for all values $\gamma>0$. The results in \cite{compositionpaper} are again stated for local classes $\mathcal{E}^{\on{loc}}$, but since we have shown there estimates for the underlying weight structures the results hold true for the ultraholomorphic class $\mathcal{A}$ as well by the analogous definition.
\qed\enddemo

Using this result and the proof of Corollary \ref{Theorem321cor} we can restate \cite[Cor. 3.4.1]{Thilliezdivision} independently by using the associated weight function.

\begin{corollary}\label{Theorem321Beurlingcor}
Let $\widehat{M}=(p!M_p)_{p\in\NN}$ be given such that $\widehat{M}\in\hyperlink{LCset}{\mathcal{LC}}$, $\gamma(\omega_{\widehat{M}})>1$ and satisfying \hyperlink{mg}{$(\on{mg})$}. Then for any $\gamma>0$ satisfying $0<\gamma<\gamma(\omega_{\widehat{M}})-1=\gamma(\omega_M)=\gamma(M)=\gamma(\widehat{M})-1$ the Borel map $\mathcal{B}:\mathcal{A}_{(\widehat{M})}(S_{\gamma})\longrightarrow\Lambda^1_{(\widehat{M})}$ is surjective.
\end{corollary}

We close this section with the following observations: In both Theorem \ref{Theorem321Beurling} and Corollary \ref{Theorem321Beurlingcor} (and also in \cite[Cor. 3.4.1]{Thilliezdivision}) using the above proofs and techniques one does not get any information about the existence of a continuous linear extension operator in the ultraholomorphic Beurling case and which should be compared in the single weight sequence situation with the results from \cite[Section 4]{Schmetsvaldivia00}.

In the weight function case, $\omega$ is said to be a $(DN)$-weight if
\begin{equation}\label{DNweight}
\forall\;C\ge 1\;\exists\delta>0\;\exists\;t_0>0\;\forall\;t\ge t_0:\;\;\omega^2(t)\le\omega(Ct)\omega(\delta t).
\end{equation}
In \cite[Corollary 3.12]{MeiseTaylor88} and in \cite{MeiseTaylor89} it was shown that in the ultradifferentiable Beurling case there does exist a continuous linear extension operator if and only if \eqref{DNweight} holds true. By \cite[Lemma 19]{BonetMeiseMelikhov07} each $\omega\in\hyperlink{omset1}{\mathcal{W}}$ with \hyperlink{omsnq}{$(\omega_{\text{snq}})$} and \hyperlink{om6}{$(\omega_6)$} has \eqref{DNweight}; more precisely we can see that for each $\omega\in\hyperlink{omset1}{\mathcal{W}}$ with \hyperlink{om6}{$(\omega_6)$} conditions \hyperlink{omsnq}{$(\omega_{\text{snq}})$} and \eqref{DNweight} are equivalent: The remaining implication follows by the fact that \eqref{DNweight} is characterizing the existence of an extension operator for $\mathcal{E}_{(\omega)}$, since $\mathcal{E}_{(\omega)}=\mathcal{E}_{(W^x)}$ holds for all $x>0$ and $W^x$ has \hyperlink{mg}{$(\on{mg})$} (see Remark \ref{importantremark} $(b)$), the characterization \cite[Theorem 3.1 (b)]{petzsche} applied to $W^x$ and finally the fact that in this situation \hyperlink{gamma1}{$(\gamma_1)$} for $W^x$ is equivalent to $\gamma(\omega_{W^x})>1$ (and so to \hyperlink{omsnq}{$(\omega_{\text{snq}})$}) for $\omega_{W^x}$, see \cite[Section 4]{firstindexpaper}. Consequently $\omega$ has \hyperlink{omsnq}{$(\omega_{\text{snq}})$} too by recalling $(vi)$ in Lemma \ref{assofuncproper}.\vspace{6pt}

Given $\widehat{M}$ as in Corollary \ref{Theorem321Beurlingcor} we can see that the associated weight function $\omega_{\widehat{M}}$ satisfies \eqref{DNweight}, in particular for any $\widehat{M}\in\hyperlink{SRset}{\mathcal{SR}}$ the function $\omega_{\widehat{M}}$ is a $(DN)$-weight. It is well-known that each $\sigma_s(t):=\max\{0,\log(t)^s\}$, $s>1$, is not a $(DN)$-weight.

\appendix

\section{A class of weight functions satisfying $\gamma(\omega)=+\infty$}\label{omega7section}
In this section we present a class of weight functions such that Theorem \ref{Theorem321} holds true for any opening $\gamma>0$ and the technical procedure in Section \ref{ramificationpreparation} becomes superfluous in the sense that we can also work in case 2 with the matrix $\Omega$ associated to $\omega$ directly. In particular all results from this section hold true for the weights $\sigma_s(t):=\max\{0,\log(t)^s\}$, $s>1$. More precisely we consider weights satisfying\vspace{6pt}

\centerline{$\omega\in\hyperlink{omset0}{\mathcal{W}_0}$ and \hyperlink{om7}{$(\omega_7)$}.}\vspace{6pt}

We point out that any nondecreasing $\omega$ with \hyperlink{om7}{$(\omega_7)$} satisfies \hyperlink{om1}{$(\omega_1)$} since, denoting by $H$ the constant arising in \hyperlink{om7}{$(\omega_7)$}, for any $t\ge 2H^2$ we have $\omega(2t)\le\omega((t/H)^2)\le C\omega(H(t/H))+C=C\omega(t)+C$.

In \cite{Franken95} for a weight function the following condition has been introduced:
\begin{equation}\label{franken}
\forall\;\gamma>1\;\exists\;L>1\;\forall\;t\ge 0:\;\;\;\omega(t^{\gamma})\le L(\omega(t)+1).
\end{equation}
For any $\omega\in\hyperlink{omset0}{\mathcal{W}_0}$, \eqref{franken} is equivalent to \hyperlink{om7}{$(\omega_7)$}: Obviously \eqref{franken} implies \hyperlink{om7}{$(\omega_7)$} by taking $\gamma=2$. Conversely given $\gamma>1$, by iterating \hyperlink{om7}{$(\omega_7)$} we have $\omega(t^{\gamma})\le\omega(t^{2^n})\le L(\omega(t)+1)$ for some $L>1$ and all $t\ge 1$ (see also the proof of $\gamma(\omega)=+\infty$ in Lemma \ref{omega7} below).\vspace{6pt}

Let $\omega\in\hyperlink{omset0}{\mathcal{W}_0}$ with \hyperlink{om7}{$(\omega_7)$} be given. By \cite[Lemma 5.9 $(5.12)$]{compositionpaper} and as we have shown in \cite[Theorem 5.14 $(4)$]{compositionpaper} for the classes $\mathcal{E}^{\text{loc}}$, we get $\mathcal{A}_{\{\omega\}}(S)=\bigcup_{x\in\RR_{>0}}\mathcal{A}_{\{W^x\}}(S)=\bigcup_{x\in\RR_{>0}}\mathcal{A}_{(W^x)}(S)$ and $\mathcal{A}_{(\omega)}(S)=\bigcap_{x\in\RR_{>0}}\mathcal{A}_{(W^x)}(S)=\bigcap_{x\in\RR_{>0}}\mathcal{A}_{\{W^x\}}(S)$ for any (unbounded) sector $S$ respectively for the class $\Lambda^n$, too.

The following result summarizes properties and consequences obtained by \hyperlink{om7}{$(\omega_7)$}, for $(i)$ see also \cite[Lemmas 3.6.1, 5.4.1]{dissertation}.

\begin{lemma}\label{omega7}
Let $\omega\in\hyperlink{omset0}{\mathcal{W}_0}$ be given, $\Omega=\{W^l: l>0\}$ the associated weight matrix. The following are equivalent:
\begin{itemize}
\item[$(i)$] $\omega$ satisfies \hyperlink{om7}{$(\omega_7)$},

\item[$(ii)$] $\exists\;A\ge 1\;\exists\;B\ge 1\;\forall\;l>0\;\exists\;C_l\ge 1\;\forall\;j\in\NN: W^l_{2j}\le C_lB^jW^{Al}_j$,

\item[$(iii)$] $\exists\;A\ge 1\;\exists\;B\ge 1\;\forall\;l>0\;\exists\;C_l\ge 1\;\forall\;j\in\NN:\;\;(W^l_j)^2\le C_lB^j W^{Al}_j$.
\end{itemize}
In case if any of the previous condition holds true we have $\gamma(\omega)=+\infty$, in particular $\omega$ is a strong weight, and moreover $\Omega\hyperlink{Mroumapprox}{\{\approx\}}\{(p!W^l_p)_{p\in\NN}: l>0\}$.
\end{lemma}

However, \hyperlink{om7}{$(\omega_7)$} is only a sufficient condition to have $\gamma(\omega)=+\infty$. One can check that the functions $\omega(t):=\exp(\beta\log(1+t)^{\alpha})=t^{\beta\log(1+t)^{\alpha-1}}$ with $\beta>0$ and $0<\alpha<1$, $t\ge t_0$ large, and which have been mentioned in \cite[Example 1.5]{BonetBraunMeiseTaylorWhitneyextension}, do satisfy $\gamma(\omega)=+\infty$ but \hyperlink{om7}{$(\omega_7)$} is violated (and \hyperlink{om6}{$(\omega_6)$} too).

\demo{Proof}
First we translate \hyperlink{om7}{$(\omega_7)$} into a property for $\varphi_{\omega}=\omega\circ\exp$: Put $t=\exp(s)$, then $t^2=\exp(2s)$ and for $h=\log(H)$ we get $Ht=\exp(\log(H))\exp(s)=\exp(h+s)$. So \hyperlink{om7}{$(\omega_7)$} means $\varphi_{\omega}(2s)\le C\varphi_{\omega}(h+s)+C$. Then apply the Legendre-Fenchel-Young-conjugate of $\varphi_{\omega}$ on both sides, i.e. for $x\ge 0$ the left-hand side gives
\begin{align*}
(\varphi_{\omega}(2\cdot))^{*}(x)&=\sup_{y\ge 0}\{xy-\varphi_{\omega}(2y)\}=\sup_{y'\ge 0}\left\{\frac{x}{2}y'-\varphi_{\omega}(y')\right\}=\varphi^{*}_{\omega}\left(\frac{x}{2}\right),
\end{align*}
where we have put $y':=2y$. For the right-hand side we get:
\begin{align*}
&(C\varphi_{\omega}(h+\cdot)+C)^{*}(x)=\sup_{y\in\RR}\{xy-C\varphi_{\omega}(h+y)-C\}=C\sup_{z\in\RR}\left\{\frac{x}{C}(z-h)-\varphi_{\omega}(z)\right\}-C
\\&
=C\sup_{z\in\RR}\left\{\frac{x}{C}z-\varphi_{\omega}(z)\right\}-hx-C=C\sup_{z\ge 0}\left\{\frac{x}{C}z-\varphi_{\omega}(z)\right\}-hx-C=C\varphi^{*}_{\omega}\left(\frac{x}{C}\right)-hx-C,
\end{align*}
where we have put $z:=h+y$. Hence \hyperlink{om7}{$(\omega_7)$} implies
\begin{equation}\label{omega7equ1}
\exists\;C,H\ge 1\;\forall\;x\ge 0:\;\;\;C\varphi^{*}_{\omega}\left(\frac{x}{C}\right)\le\varphi^{*}_{\omega}\left(\frac{x}{2}\right)+\log(H)x+C.
\end{equation}

$(i)\Longrightarrow(ii)$ In \eqref{omega7equ1} we put $x=2lj$ for $l>0$ and $j\in\NN$. Then divide by $l$ and apply $\exp$ to get by definition
$$\exists\;C,H\ge 1:\;\forall\;l>0\;\forall\;j\in\NN: W^{l/C}_{2j}\le\exp\left(\frac{C}{l}\right)(H^2)^jW^l_j,$$
so put $A:=C$, $B:=H^2$ and $C_l:=\exp(C/l)$.\vspace{6pt}

$(ii)\Longrightarrow(i)$ After applying $\log$ we obtain for all $j\in\NN$:
$$\frac{1}{l}\varphi^{*}_{\omega}(2lj)\le\log(C_l)+j\log(B)+\frac{1}{Al}\varphi^{*}_{\omega}(Alj).$$
Replace $j\in\NN$ by $y\ge 0$ and apply the Legendre-Fenchel-Young-conjugate to both sides for $x\ge 0$. The left-hand side gives
\begin{align*}
\left(\frac{1}{l}\varphi^{*}_{\omega}(2l\cdot)\right)^{*}(x)&=\sup_{y\ge 0}\left\{xy-\frac{1}{l}\varphi^{*}_{\omega}(2ly)\right\}=\frac{1}{l}\sup_{z\ge 0}\left\{\frac{lxz}{2l}-\varphi^{*}_{\omega}(z)\right\}=\frac{1}{l}\varphi_{\omega}^{**}\left(\frac{x}{2}\right)=\frac{1}{l}\omega(\sqrt{\exp(x)}),
\end{align*}
whereas the right-hand side gives
\begin{align*}
&\left(\log(C_l)+\cdot\log(B)+\frac{1}{Al}\varphi^{*}_{\omega}(Al\cdot)\right)^{*}(x)=\sup_{y\ge 0}\left\{xy-y\log(B)-\frac{1}{Al}\varphi^{*}_{\omega}(Aly)\right\}-\log(C_l)
\\&
=\frac{1}{Al}\sup_{z\ge 0}\left\{\frac{z}{Al}(x-\log(B))Al-\varphi^{*}_{\omega}(z)\right\}-\log(C_l)=\frac{1}{Al}\varphi^{**}_{\omega}(x-\log(B))-\log(C_l)
\\&
=\frac{1}{Al}\omega\left(\frac{\exp(x)}{B}\right)-\log(C_l).
\end{align*}
We summarize:
\begin{align*}
&\frac{1}{l}\omega(\sqrt{\exp(x)})=\sup_{y\ge 0}\left\{xy-\frac{1}{l}\varphi^{*}_{\omega}(2ly)\right\}\ge\sup_{j\in\NN}\left\{xj-\frac{1}{l}\varphi^{*}_{\omega}(2lj)\right\}
\\&
\ge\sup_{j\in\NN}\left\{xj-j\log(B)-\frac{1}{Al}\varphi^{*}_{\omega}(Alj)\right\}-\log(C_l)\underbrace{\ge}_{(\circ)}\frac{1}{2}\sup_{y\ge 0}\left\{xy-y\log(B)-\frac{1}{Al}\varphi^{*}_{\omega}(Aly)\right\}-\log(C_l)
\\&
=\frac{1}{2Al}\omega\left(\frac{\exp(x)}{B}\right)-\log(C_l).
\end{align*}
$(\circ)$ holds by \cite[Theorem 4.0.3, Lemma 5.1.3]{dissertation} respectively \cite[Lemma 5.7]{compositionpaper}. Put $n:=2Al$, then $\omega_n:=\frac{1}{n}\omega\hyperlink{sim}{\sim}\omega_{W^n}$ and $\omega_n(t)\le 2\omega_{W^n}(t)$ for all $t\ge\frac{W^n_2}{W^n_1}$. Put $t:=\exp(x)$ and so for any $\frac{t}{B}\ge\frac{W^n_2}{W^n_1}$ we have shown $\omega(\frac{t}{B})\le\frac{n}{l}\omega(\sqrt{t})+Al\log(C_l)=2A\omega(\sqrt{t})+Al\log(C_l)$.\vspace{6pt}

$(ii)\Longrightarrow(iii)$ By using log-convexity for $W^l$ we get $(W^l_j)^2\le W^l_{2j}\le C_lB^jW^{Al}_j$.

$(iii)\Longrightarrow(ii)$ By \eqref{newmoderategrowth} we have $W^l_{2j}\le(W^{2l}_j)^2\le C_lB^j W^{A2l}_j$.\vspace{6pt}


To show $\gamma(\omega)=+\infty$, let $\gamma>1$ be given, arbitrary but from now on fixed, and $C,H\ge 1$ shall denote the constants coming from \hyperlink{om7}{$(\omega_7)$}. Iterating \hyperlink{om1}{$(\omega_1)$} $n$-times, where $n\in\NN_{>0}$ is chosen minimal such that $H\le 2^n$ holds, we get
$$\exists\;D\ge 1\;\forall\;t\ge 0:\;\;\;\omega(Ht)\le\omega(2^nt)\le D\omega(t)+D,$$
hence
$$\exists\;D_1\ge 1\;\forall\;t\ge 0:\;\;\;\omega(t^2)\le CD\omega(t)+CD+C\le D_1\omega(t)+D_1,$$
which means $\limsup_{t\rightarrow+\infty}\frac{\omega(t^2)}{\omega(t)}\le D_1$ (put $D_1:=CD+C\ge 1$). Then defining $K:=D_1+1>1$ and for all $t\ge 0$ sufficiently large depending on given $\gamma$ we have $\omega(K^{\gamma}t)\le\omega(t^2)$ because $\omega$ is nondecreasing. Consequently
$$\limsup_{t\rightarrow+\infty}\frac{\omega(K^{\gamma}t)}{\omega(t)}\le\limsup_{t\rightarrow+\infty}\frac{\omega(t^2)}{\omega(t)}\le D_1<K.$$

Finally, let us prove that $\Omega\hyperlink{Mroumapprox}{\{\approx\}}\{(p!W^l_p)_{p\in\NN}: l>0\}$. Since $\omega$ is a strong weight it satisfies \hyperlink{om5}{$(\omega_5)$}, see Proposition \ref{Prop13MT88}. So by \cite[Proposition 4.6, Corollary 5.15]{compositionpaper} for all $l>0$ there exists some $D_l\ge 1$ such that $j!\le D_l^jW^l_j$ holds for all $j\in\NN$ (apply Stirling's formula). Then applying $(iii)$ we get $j!W^l_j\le D_l^j(W^l_j)^2\le C_l(D_lB)^jW^{Al}_j$ for all $j\in\NN$ and constants $B,A\ge 1$ both not depending on given $l$. Note that without loss of generality we can take $D_l=D_1$ for any $l\ge 1$ since in this case $W^l\ge W^1$.
\qed\enddemo

\begin{remark}\label{remarkomega7}
\begin{itemize}
\item[$(a)$] As already pointed out in \cite[Lemma 5.9]{compositionpaper} any $\omega\in\hyperlink{omset}{\mathcal{W}}$ with \hyperlink{om7}{$(\omega_7)$} can never satisfy \hyperlink{om6}{$(\omega_6)$} and consequently $\mathcal{A}_{\{\omega\}}=\mathcal{A}_{\{M\}}$ can not hold for any single weight sequence $M\in\hyperlink{LCset}{\mathcal{LC}}$.

\item[$(b)$] By iterating $(iii)$ we see that this property is equivalent to the fact that $\Omega\hyperlink{Mroumapprox}{\{\approx\}}\{(W^l)^q: l>0\}$ for any $q>0$. Hence for any $\omega\in\hyperlink{omset}{\mathcal{W}}$ with \hyperlink{om7}{$(\omega_7)$} we see that $\Omega\hyperlink{Mroumapprox}{\{\approx\}}\{(p!(w^l_p)^q)_{p\in\NN}: l>0\}$ for any $q>0$, consequently Section \ref{ramificationpreparation} (necessary to prove case 2 in Theorem \ref{Theorem321}) becomes superfluous since we can use the same weight matrix and work with
    $$\Omega\hyperlink{Mroumapprox}{\{\approx\}}\{(p!(w^l_p)^q)_{p\in\NN}: l>0\}\hyperlink{Mroumapprox}{\{\approx\}}\widehat{\mathcal{S}}^{q}$$ directly.

\item[$(c)$] We also get in this situation that each $\omega_{w^x}$ is a strong weight and $\omega_{W^x}\hyperlink{sim}{\sim}\omega\hyperlink{sim}{\sim}\omega_{w^x}$ for all $x>0$ (and which also implies $\gamma(\omega)=+\infty$ by Corollary \ref{growthindexgammacorollary}). Note that $\omega_{M^x}(t)\le\omega_{m^x}$ holds for any $x>0$, hence $\omega_{w^x}\hyperlink{ompreceq}{\preceq}\omega_{W^x}\hyperlink{sim}{\sim}\omega$. On the other hand, given $x>0$, we get $j!W^x_j\le B^jW^{Ax}_j\Leftrightarrow W^x_j\le B^jw^{Ax}_j$ for some $A,B\ge 1$ and all $j\in\NN$, where $A$ is not depending on $x$. This implies $\omega_{w^{Ax}}(t)\le\omega_{W^x}(Bt)$ for all $t\ge 0$. Since $\omega_{W^x}\hyperlink{sim}{\sim}\omega$ we get \hyperlink{om1}{$(\omega_1)$} for $\omega_{W^x}$ and so $\omega_{w^{Ax}}(t)\le B'\omega_{W^x}(t)+B'$ for some $B'\ge 1$ and all $t\ge 0$. Consequently $\omega\hyperlink{sim}{\sim}\omega_{W^{x/A}}\hyperlink{ompreceq}{\preceq}\omega_{w^x}$ holds true.

\item[$(d)$] Finally in this case we can replace $\Omega$ by an equivalent matrix $\mathcal{M}$ consisting of only strongly nonquasianalytic weight sequences, i.e. \hyperlink{beta1}{$(\beta_1)$} holds true for each $M^x\in\mathcal{M}$. This follows by having $\Omega\hyperlink{Mroumapprox}{\{\approx\}}\{(p!^2W^l_p)_{p\in\NN}: l>0\}$ and by taking into account the arguments given in \cite[Example 7.10 $(ii)$]{sectorialextensions}.
\end{itemize}
\end{remark}

\textbf{Acknowledgements}: The first two authors are partially supported by the Spanish Ministry of Economy, Industry and Competitiveness under the project MTM2016-77642-C2-1-P. The first author is partially supported by the University of Valladolid through a Predoctoral Fellowship (2013 call) co-sponsored by the Banco de Santander. The third author is supported by FWF-Project J~3948-N35, as a part of which he is an external researcher at the Universidad de Valladolid (Spain) for the period October 2016 - September 2018.\par

\bibliographystyle{plain}
\bibliography{Bibliography}

\vskip1cm

\textbf{Affiliation}:\\
J.~Jim\'{e}nez-Garrido, J.~Sanz:\\
Departamento de \'Algebra, An\'alisis Matem\'atico, Geometr{\'\i}a y Topolog{\'\i}a, Universidad de Valladolid\\
Facultad de Ciencias, Paseo de Bel\'en 7, 47011 Valladolid, Spain.\\
Instituto de Investigaci\'on en Matem\'aticas IMUVA\\
E-mails: jjjimenez@am.uva.es (J.~Jim\'{e}nez-Garrido), jsanzg@am.uva.es (J. Sanz).
\\
\vskip.5cm
G.~Schindl:\\
Departamento de \'Algebra, An\'alisis Matem\'atico, Geometr{\'\i}a y Topolog{\'\i}a, Universidad de Valladolid\\
Facultad de Ciencias, Paseo de Bel\'en 7, 47011 Valladolid, Spain.\\
E-mail: gerhard.schindl@univie.ac.at.

\end{document}